\documentclass[10pt]{amsart}
\usepackage{amssymb}
\usepackage{epsfig}
\usepackage{amsfonts}
\usepackage{amscd}

\newcommand{\mbotreix}{\overline{\mathcal{M}}_{0,3}(X,d)}
\newcommand{\mbotreig}{\overline{\mathcal{M}}_{0,3}(Gr(p,m),d)}

\newcommand{\z}{\mathbb{Z}}

\newcommand{\cx}{\mathbb{C}}

\newtheorem{thm}{\textbf{Theorem }}
\newtheorem{prop}{\textbf{Proposition }}[section]
\newtheorem{cor}[prop]{\textbf{Corollary}}
\newtheorem{lemma}[prop]{\textbf{Lemma}}

\newcommand{\eqqh}{{QH}^{*}_{T}(Gr(p,m))}
\newcommand{\eqqhx}{{QH}^{*}_{T}(X)}

\begin{document}
\title{Equivariant quantum Schubert Calculus}
\author{Leonardo Constantin Mihalcea }
\date{\today} \maketitle
\begin{abstract}
We study the $T-$equivariant quantum cohomology of the
Grassmannian. We prove the vanishing of a certain class of
equivariant quantum Littlewood-Richardson coefficients, which
implies an equivariant quantum Pieri rule. As in the equivariant
case, this implies an algorithm to compute the equivariant quantum
Littlewood-Richardson coefficients.
\end{abstract}
\section{Introduction}

A deformation of the integral cohomology of $X=Gr(p,m)$, the
Grassmannian of $p-$planes in $\cx^m$, has been recently
constructed. It is the (small) quantum cohomology of $X$, which is
a graded algebra over $\z[q]$, where the (complex) degree of $q$
is $m$. It has a $\z[q]-$basis consisting of Schubert classes $\{
\sigma_\lambda \}$ indexed by partitions
$\lambda=(\lambda_1,...,\lambda_p)$ included in the $p \times
(m-p)$ rectangle (i.e. $m-p \geqslant \lambda_1 \geqslant
\lambda_2 \geqslant ... \geqslant \lambda_p \geqslant 0$). The
multiplication is governed by the quantum Littlewood-Richardson
coefficients, a special case of the 3-point Gromov-Witten
invariants, which encode enumerative properties of the variety $X$
(\cite{W,KM,FP,Be,C,BCF,FGP,Bu1,Po}).
\\ \indent The purpose of this paper is to study a $T-$equivariant
version of the quantum cohomology, where $T \simeq (\cx^\star)^m$
is the torus of diagonal $m \times m$ invertible matrices acting
on $X$. This was introduced by Givental and Kim (\cite{GK}) to
study the (ordinary) quantum cohomology algebra. It is named the
equivariant quantum cohomology and it is a deformation of both
quantum and $T-$equivariant cohomology of $X$. It has a structure
of a graded $\Lambda[q]$ algebra, where $\Lambda$ denotes the
$T-$equivariant cohomology of a point, which is identified with
the polynomial ring $\z[T_1,...,T_m]$ (for a geometric description
of $T_i$ see \S 2.2 below). Additively, it has a
$\Lambda[q]$-basis $\{ \sigma_\lambda \}$ indexed by partitions
$\lambda$ included in the $p \times (m-p)$ rectangle. The
multiplication, denoted $\circ$, is determined by the
\textit{equivariant quantum Littlewood-Richardson coefficients}
(EQLR) $c_{\lambda\mu}^{\nu,d}$, a special case of the 3-point
\textit{equivariant Gromov-Witten invariants}, introduced in
\cite{GK} (for more general varieties). Explicitly:
\[\sigma_\lambda \circ \sigma_\mu = \sum_{d \geqslant 0}\sum_{\nu}
c_{\lambda,\mu}^{\nu,d}q^d \cdot \sigma_{\nu} \] \indent By
definition, the EQLR coefficient $c_{\lambda,\mu}^{\nu,d}$ is a
homogeneous polynomial in $\Lambda$ of degree $|\lambda|+|\mu| -
|\nu|-md$, where $|\lambda|=\lambda_1+...+\lambda_p$ denotes the
weight of the partition $\lambda$. If $c_{\lambda,\mu}^{\nu,d}$
has polynomial degree $0$, it is equal to the quantum LR
coefficient $c_{\lambda\mu}^{\nu,d}$, while if $d=0$ the EQLR
coefficient is equal to the equivariant LR coefficient
$c_{\lambda\mu}^\nu$ (\cite{Kim1}). There are explicit formulae
for the equivariant (\cite{KT1,MS}, Prop. \ref{equivariant
reccurence} below) and quantum Littlewood-Richardson coefficients
(\cite{BCF,BKT1}). Therefore we are mainly interested in ``mixed"
EQLR coefficients, i.e. those with positive polynomial degree for
which $d>0$.

\indent Our goal is to give an equivariant quantum Pieri rule for
the equivariant quantum cohomology of $X$ and an effective
algorithm to compute the EQLR coefficients. The key to that is to
employ Buch's notions of \textit{span} and \textit{kernel} of a
stable map to $X$ (\cite{Bu1}) as well as classical Schubert
calculus to obtain vanishing
properties for some of these coefficients.\\

\noindent \textbf{Conventions:} 1. Unless otherwise specified, all
the partitions used in this paper are included in the $p \times
(m-p)$ rectangle.

2. The {\it degree} of $c_{\lambda,\mu}^{\nu,d}$ is $d$. \\

\noindent \textbf{Notations:} Write $\alpha \rightarrow \beta$ if $\beta$
is included in $\alpha$ and the Young diagram of $\alpha$ has one more
box than the diagram of $\beta$. Denote by $\alpha^-$ (resp. $\alpha^+$)
the partition obtained from $\alpha$ by removing (resp. by adding) $m-1$
boxes from (resp. to) its border rim (recall that the border rim of a
Young diagram is the set of boxes that intersect the diagram's SE
border). The partition $\alpha^+$ can also be defined by the equivalent
property $(\alpha^+)^-=\alpha$. If $\alpha=(\alpha_1,...,\alpha_p)$, note
that $\alpha^-$ exists only if $\alpha_1=m-p$ and $\alpha_p>0$, and that
$\alpha^+$ is included in the $p \times (m-p)$ rectangle only if
$\alpha_1<m$ and $\alpha_p=0$.

Write $\alpha^\vee$ for the partition dual to $\alpha$, i.e. the
partition whose Young diagram is the rotation with 180 degrees of
the complement of the Young diagram of $\alpha$, in the given
rectangle. $\alpha'$ will denote the partition conjugate to
$\alpha$, i.e. the partition in the $(m-p) \times p$ rectangle
whose Young diagram is the transpose of the Young diagram of
$\alpha$.

The zero partition (i.e. the partition with all parts of length
zero) is denoted by $(0)$. Note that $\sigma_{(0)}$ is the unit in
the cohomology of the Grassmannian. The partition $(1,0,...,0)$ is
denoted by $\Box$.

\subsection{Statement of the results}
The main result of this paper is an equivariant quantum
Pieri\begin{footnote} [1] {In equivariant Schubert calculus the
special multiplication by $\sigma_\Box$ has a more prominent role
than in the classical one, therefore it is referred to as the
``Pieri rule". A more general equivariant Pieri rule, involving
multiplication by $\sigma_{(k)}$ or $\sigma_{(1)^k}$ has been
recently obtained in \cite{R} for the complete flag manifold.}
\end{footnote}rule for the Grassmannians. The rule must be a
deformation of the corresponding equivariant and quantum rules
(\cite{Be,KT1}, (see Sections 2.2 and 2.3 for details) and it is
the simplest possible under these circumstances.
There are no ``mixed" terms:   \\

\noindent \textbf{Theorem (Equivariant Quantum Pieri rule)}
\textit{ The equivariant quantum multiplication for the
Grassmannian $Gr(p.m)$ satisfies the following formula:
 \[\sigma_\lambda \circ \sigma_\Box =
\sum_{\mu \rightarrow \lambda} \sigma_{\mu} + c_{\lambda,\Box
}^\lambda \sigma_\lambda + q \sigma_{\lambda^-} \] where
$c_{\lambda,\Box}^\lambda$ is the equivariant LR coefficient,
known by Prop. \ref{equivariant reccurence} below. The last term
is omitted if
$\lambda^-$ does not exist.}\\

\noindent \textit{Remark:} Corollary 7 in \cite{KiMa} gives
another formula which is referred therein as the ``equivariant
quantum Pieri rule". It deals with a Pieri multiplication in the
presentation of the $U_m$-equivariant quantum cohomology algebra
of the complete flag manifold $Fl(m)$ obtained in \cite{GK}. Since
the (equivariant) quantum cohomology is not functorial, this
formula does not imply one for the Grassmannian.\\

The theorem follows from a key vanishing condition of the EQLR coefficients: \\

\noindent \textbf{Main Lemma} \textit{Let $\lambda,\mu,\nu$ be
three partitions included in $p \times (m-p)$ rectangle and let
$d$ be a positive integer. Suppose that $|\lambda|+d^2 >
|\nu|+md$. Then
$c_{\lambda,\mu}^{\nu,d}=0$.}\\

In particular, $c_{\lambda,\mu}^{\nu,d}=0$ if $|\lambda|+|\mu|>
|\nu|+md $, and $\mu$ is included in the $d \times d$ square. This
applies to the mixed EQLR coefficients $c_{\lambda,\Box}^{\nu,d}$,
implying the EQ Pieri rule.

An algebraic consequence of the equivariant quantum Pieri rule and
of an associativity equation of the equivariant quantum cohomology
is a recursive formula satisfied by the EQLR coefficients. It
expresses $c_{\lambda,\mu}^{\nu,d}$ as a combination of EQLR
coefficients with degree $d-1$, and EQLR coefficients with the
same degree $d$, but with polynomial degree one larger. The
formula is a generalization of a recursive formula for the
equivariant LR coefficients (see \cite{MS,O,KT1}, or Prop.
\ref{equivariant
reccurence} below). \\

\noindent \textbf{Corollary} \textit{The EQLR coefficients satisfy
the following formula:\[(c_{\nu,\Box }^\nu - c_{\lambda,\Box
}^\lambda)\cdot c_{\lambda,\mu}^{\nu,d} = \sum_{\delta \rightarrow
\lambda} c_{\delta,\mu}^{\nu,d} - \sum_{\nu \rightarrow \zeta}
c_{\lambda,\mu}^{\zeta,d} + c_{\lambda^-, \mu}^{\nu,d-1} -
c_{\lambda,\mu}^{\nu^+,d-1} \] for any partitions
$\lambda,\mu,\nu$ and any nonnegative integer $d$, where
$c_{\alpha,\Box}^\alpha$ is the equivariant LR coefficient given
in Prop. \ref{equivariant reccurence}. The third (resp. the
fourth) term in the right side is omitted if $\lambda^-$ (resp.
$\nu^+$) does not exist in the $p \times (m-p)$ rectangle. Both
these terms are omitted if $d=0$.} \\

This formula is the main ingredient in the proof of an algorithm,
which shows that the EQLR coefficients are determined by their
usual commutativity equation, by those appearing in the
multiplication with the unit $\sigma_{(0)}$, by the Pieri
coefficients, and by the formula from
the previous Corollary (for the precise statement, see Thm. \ref{recalg} in \S 7 below).\\

\noindent \textit{Acknowledgements:} I would like to thank to my
advisor, Prof. W. Fulton, for pointing me to this research area
and for his patient guidance, which greatly improved the
presentation of this paper. This is part of my thesis.
\section{Preliminaries}
In this section we recall some basic facts about the classical,
equivariant and quantum cohomology of the Grassmannian, which are
needed later in the paper.
\subsection{Classical cohomology of the Grassmannian}
Let $Gr(p,m)$ be the Grassmannian of $p-$planes in $\cx^m$. Fix a
complete flag $F_\bullet = 0 \subset F_1 \subset ... \subset F_m=
\cx^m$. Let $\lambda$ be a partition included in the $p \times
(m-p)$ rectangle and $\Omega_\lambda(F_\bullet)$ be the Schubert
variety determined by $F_\bullet$ and
$\lambda=(\lambda_1,...,\lambda_p)$ i.e.
\[ \Omega_\lambda(F_\bullet)= \{ V \in Gr(p,m): \dim (V \cap
F_{m-p+i-\lambda_i}) \geqslant i \} \] Denote by $\sigma_\lambda$
the cohomology class in $H^{2|\lambda|}(Gr(p,m))$ determined by
$\Omega_\lambda(F_\bullet)$. It is well-known that the classes
$\sigma_\lambda$ do not depend on the choice of the flag
$F_\bullet$, and that they form a $\z-$basis for the integral
cohomology of $Gr(p,m)$ (see \cite{F1} Part III for an exposition
about the subject). The multiplication in the cohomology ring is
determined by the Littlewood-Richardson (LR) coefficients
$c_{\lambda,\mu}^\nu$, which are positive integers, counting the
number of points in the intersection of the Schubert varieties
$\Omega_\lambda(F_\bullet)$, $\Omega_\mu(G_\bullet) $ and $
\Omega_{\nu^\vee}(H_\bullet)$, where $F_\bullet$,$G_\bullet$ and
$H_\bullet$ are three general flags and $\nu^\vee$ is the dual
partition of $\nu$. The coefficients are 0 if $|\lambda| + |\mu|
\neq |\nu|$. Geometrically,
\[ c_{\lambda,\mu}^{\nu} = \pi_\star (\sigma_\lambda \cup \sigma_\mu
\cup \sigma_{\nu^\vee})\] where $\pi_\star:H^i(Gr(p,m)) \to
H^{i-2p(m-p)}(pt)$ is the Gysin map associated to the structure
morphism $\pi:Gr(p,m) \longrightarrow pt$ (see Appendix for
details). Positive combinatorial formulae for these coefficients,
found in the literature as Littlewood-Richardson rules, are known
(see e.g. \cite{F1} Part I or \cite{S} and references therein).

\subsection{Equivariant cohomology} Let $T$ be the $m-$dimensional
complex torus $(\cx^\star)^m$ and $X$ an algebraic variety with a
$T-$action. Let $p : ET \longrightarrow BT$ be the universal
$T-$bundle. There is an induced $T-$action on the product $ET
\times X$ given by $t\cdot (e,x)=(t^{-1}e,tx)$. This determines a
quotient space $X_T=ET \times_T X$, the homotopic quotient of $X$.
It is an $X-$bundle over $BT$. The $T-$equivariant cohomology of
$X$, denoted by $H^\star_T (X)$, is by definition the ordinary
(integral) cohomology of $X_T$. The $X-$bundle projection $\pi:
X_T \longrightarrow BT$ gives $H^\star_T (X)$ a structure of
$\Lambda-$algebra, where $\Lambda$ denotes
$H^\star_T(pt)=H^\star(BT)$ (for details see
\cite{AB,GKM,Br1,Br}).

The topological spaces $ET$ and $BT$ are infinite dimensional, so
in particular they are not algebraic varieties. Nevertheless, one
can consider the direct system of finite-dimensional $T-$bundles
$p: ET_n \longrightarrow BT_n$ given by \[ \prod_{i=1}^m (\cx^n
\smallsetminus \{0\}) \longrightarrow \prod_{i=1}^{m}
\mathbb{P}^{n-1}\] The ordering on the bundles is given by
inclusion (for the construction above or similar ones see e.g.
\cite{Br}, \cite{H}, Ch.4 \S 11, or Ch. 7, \cite{EG} \S 3.1). Let
$X_{T,n}:=ET_n \times_T X$ be the induced finite dimensional
approximations of $X_T$. Then one can show that $ H^i_T(X)$ is
equal to $H^i(X_{T,n})$ for $n$ large (\cite{Br}).

In particular, the equivariant cohomology of a point is
$H^i_T(pt)=H^i(BT_n)=H^i(\prod_{i=1}^m\mathbb{P}^{n-1})$ for $n$
large. By letting $n$ to go to infinity, we get that
\[ \Lambda=H^\star_T(pt)=\z[T_1,...,T_m]\] where $T_i$ has complex degree 1, and
is equal to the first Chern class
$c_1(\mathcal{O}_{\mathbb{P}^{n-1}_{(i)}}(1))$ of the line bundle
$\mathcal{O}(1)$ on $\mathbb{P}^{n-1}_{(i)}$, the $i-$th component
of the product $\prod_{i=1}^m\mathbb{P}^{n-1}$ (for a more
intrinsic definition of $H^\star_T(pt)$, see e.g. \cite{Br1,Gr}).

Equivariant cohomology has functorial properties similar to those
of ordinary cohomology. If $f: X \longrightarrow Y$ is a
$T-$equivariant map of topological spaces, it induces an pull-back
map in cohomology $f^\star:H^\star_T(Y) \longrightarrow
H^\star_T(X)$. In certain situations, for such a $T-$equivariant
map, there is also a Gysin map in cohomology
\[f_\star^T: H^i_T(X) \longrightarrow H^{i - 2d}_T(Y)\] where $d=
\dim (X) - \dim (Y)$. For the purpose of this paper, we only consider the
situation when $X$ and $Y$ are projective algebraic varieties and $Y$ is
smooth.\begin{footnote}[2]{For more general situations, such as $X$ or
$Y$ noncompact, or being able to find an ``orientation" for the map $f$,
one can consult e.g. \cite{FM}, or \cite{F3}, Ch. 19.}
\end{footnote} The definition and some of the properties of the Gysin
maps can be found in the Appendix. A more detailed treatment is
given in my thesis (\cite{Mi}).

A particular case of such an equivariant Gysin map, which will
play an important role in what follows, is the ``integration along
the fibres" $\pi_\star^T:H^i_T(X) \longrightarrow H^{i - 2 \dim
(X)}_T(pt)$ induced by the $T-$equivariant map $\pi: X
\longrightarrow pt$. It determines a $\Lambda-$pairing \[ \langle
\cdot, \cdot \rangle_T : H^\star_T(X) \otimes_\Lambda H^\star_T(X)
\longrightarrow \Lambda
\] defined by \[ \langle a,b \rangle_T = \pi_\star^T(a \cup b) \]
More about this pairing will be given in Prop. \ref{duality}.

\subsubsection{Equivariant Schubert calculus} \indent Let $X$ be the
Grassmannian $Gr(p,m)$. It inherits a diagonal $T-$action by
restriction from its $Gl(m)-$ action. Given the finite-dimensional
approximation $ET_n$ chosen above, one can show that the
$X-$bundle $\pi: X_{T,n} \longrightarrow BT_n$ is the Grassmann
bundle
\[ \mathbb{G}(p,\mathcal{O}_{(1)}(-1) \oplus ... \oplus
\mathcal{O}_{(m)}(-1)) \longrightarrow \mathbb{P}^{n-1}_{(1)}
\times ... \times \mathbb{P}^{n-1}_{(m)}\]where $\mathbb{G}(p,E)$
is the Grassmann bundle of rank $p$ subbundles of the vector
bundle $E$ (see e.g. \cite{EG}, Section 3.3, for the case
$T=\cx^\star$, $p=1$ ). In particular, note that $X_{T,n}$ is a
smooth projective variety for any positive integer $n$.

Denote by $\Omega_\lambda$ the Schubert variety of $X$ determined
by the standard flag \[ F_\bullet: 0 \subset \langle e_1 \rangle
\subset \langle e_1,e_2\rangle \subset ... \subset \langle
e_1,...,e_m \rangle = \cx^m \] and by $\widetilde{\Omega}_\lambda$
the Schubert variety determined by the opposite flag \[
F_\bullet^{opp}: 0 \subset \langle e_m \rangle \subset \langle
e_m,e_{m-1}\rangle \subset ... \subset \langle e_m,...,e_1 \rangle
= \cx^m \] These are $T-$stable varieties, so they determine
varieties $\Omega_\lambda \times_T ET_n$ in $X_{T,n}$, denoted by
$\Omega_{\lambda,n}$, and $\widetilde{\Omega}_\lambda \times_T
ET_n$, denoted by $\widetilde{\Omega}_{\lambda,n}$. These
varieties determine cohomology classes denoted by
$\sigma_{\lambda,n}^T $ respectively by
$\widetilde{\sigma}_{\lambda,n}^T$ in $H^{2|\lambda|}(X_{T,n})$.
Since these varieties are compatible as $n$ varies, they determine
equivariant cohomology classes denoted by $\sigma_\lambda^T$ and
$\widetilde{\sigma}_\lambda^T$ in $H^{2|\lambda|}_T(X)$. Note that
cohomology classes determined by the Schubert varieties depend on
the flag used to define them.

Since the classes $\sigma_\lambda$ form a $\z-$basis for the
cohomology of the fibers of $X_{T,n} \longrightarrow BT_n$, the
Leray-Hirsch theorem implies that the classes
$\sigma_{\lambda,n}^T$ form a $H^\star(BT_n)$-basis for
$H^\star(X_{T,n})$ (see \cite{H}, Ch. 16). It follows that
$\{\sigma_\lambda^T\}$ form a $\Lambda-$basis of $H^\star_T(X)$.

The structure constants of the equivariant cohomology with respect
to this basis, denoted by $c_{\lambda,\mu}^\nu$, are called the
\textit{equivariant Littlewood-Richardson coefficients}. They
agree with the classical ones when $|\lambda|+|\mu|=|\nu|$ and
they are defined by the following formulae in $H^\star_T(X)$: \[
\sigma_\lambda^T \cdot \sigma_\mu^T = \sum_{\nu}
c_{\lambda,\mu}^\nu \sigma_\nu^T
\] From definition it follows that
$c_{\lambda,\mu}^\nu=c_{\mu,\lambda}^\nu$ and that this is a homogeneous
polynomial of (complex) degree $|\lambda| + |\mu| - |\nu|$ in
$\Lambda=\z[T_1,...,T_m]$.\begin{footnote}[3] {A deep result, conjectured
by Peterson and proved by Graham (\cite{Gr}) implies that
$c_{\lambda,\mu}^\nu$ is a polynomial in $\z[T_1-T_2,...,T_{m-1}-T_m]$
with \textit{nonnegative} coefficients. A positive combinatorial formula
in this sense was then obtained in \cite{KT1}. The coefficients
$c_{\lambda,\mu}^{\nu}$ are expressed in terms of sums of weighted
puzzles, where the weight of each puzzle is a monomial in the variables
$T_1-T_2,...,T_{m-1}-T_m$ with coefficient equal to $1$.}\end{footnote}

One way to effectively compute these coefficients is a recurrence
formula, which appears in \cite{MS,O,KT1}. To state it, more
notations are needed.

Each partition $\lambda$ in the $p \times (m-p)$ rectangle is
traced out by a path starting from the NE corner of the $p \times
(m-p)$ rectangle and ending on the SW corner of the rectangle.
Define the sets $I(\lambda)$ and $J(\lambda)$ encoding the
positions of the vertical and horizontal steps of this path:
\begin{center} $ I(\lambda)=
\{ i : \textrm{ the } i-\textrm{th step of the path of } \lambda
\textrm{ is vertical } \} $ \end{center}
\begin{center} $ J(\lambda)= \{ j : \textrm{ the } j-\textrm{th step
of the path of } \lambda \textrm{ is horizontal } \} $
\end{center} For example, for $p=2,m=4$ and $\lambda=(1,1)$, the set
$I(\lambda)$ is $\{2,3\}$ while $J(\lambda)$ is $\{1,4\}$. Recall
that $\Box$ denotes the partition $(1,0,...,0)$.
\begin{prop} [\cite{MS,O,KT1}]\label{equivariant reccurence}
The equivariant LR coefficients $c_{\lambda,\mu}^\nu$ are
determined (algorithmically) by the following formulae:

(a) $ c_{\lambda,\Box }^\lambda = \sum_{i \in I(\lambda)} T_i -
\sum_{j=m-p+1}^m T_j $

(b) $c_{\lambda,\lambda}^\lambda =\prod_{i \in I(\lambda), j \in
J(\lambda), i<j} (T_i - T_j)$

(c) $(c_{\lambda,\Box }^\lambda - c_{\mu,\Box }^\mu)\cdot
c_{\lambda,\mu}^\lambda = \sum_{\delta \rightarrow \mu}
c_{\lambda,\delta}^\lambda$ for any $\lambda, \mu$ such that
$\lambda \neq \mu$.

(d) $(c_{\nu,\Box }^\nu - c_{\lambda,\Box }^\lambda)\cdot
c_{\lambda,\mu}^\nu = \sum_{\delta \rightarrow \lambda}
c_{\delta,\mu}^\nu - \sum_{\nu \rightarrow \zeta}
c_{\lambda,\mu}^{\zeta}$ for any $\lambda,\mu,\nu$ such that
$\lambda \neq \nu$.
\end{prop}
\noindent Except for parts (a) and (b), the proposition follows
immediately from the equivariant Pieri rule:
\[ \sigma_\lambda^T \cdot \sigma_\Box^T = \sum_{\delta
\rightarrow \lambda} \sigma_{\delta}^T +  c_{\lambda,\Box
}^\lambda \sigma_\lambda^T \]

There is also a geometric definition of $c_{\lambda,\mu}^\nu$. It
is the direct generalization of the geometric formula for the LR
coefficients presented in \S 2.1. As in the classical case, it
follows from a Duality theorem, stated with respect to the {\it
equivariant} Poincar\'e pairing defined in the previous section:
\begin{prop}[equivariant Duality Theorem]\label{duality} The following
formula holds in $H^0_T(pt)$:
 \[ \langle \sigma_\lambda^T ,
\widetilde{\sigma}_\mu^T \rangle_T = \pi_\star^T(\sigma_\lambda^T
\cup \widetilde{\sigma}_\mu^T)= \delta_{\lambda^\vee\mu} \] where
$\delta_{\lambda^\vee\mu}$ is the Kronecker symbol, equal to $1$
if $\lambda^\vee$ is equal to $\mu$ and $0$ otherwise.
\end{prop}
\begin{proof} This is found in the proof of Lemma
4.2, \cite{Gr} for the variety of complete flags. The result for
Grassmannians follows by pulling back using the ($T-$invariant)
projection $ pr:Fl(m) \longrightarrow Gr(p,m) $ which induces an
injective homomorphism \[pr^\star_T: H^\star_T(Gr(p,m))
\longrightarrow H^\star_T(Fl(m))\] in the equivariant cohomology.
For a direct proof see \cite{Mi}. \end{proof}

The proposition implies that \[ c_{\lambda,\mu}^\nu =
\pi_\star^T(\sigma_{\lambda}^T \cup\sigma_{\mu}^T \cup
\widetilde{\sigma}_{\nu^\vee}^T) \] This is the form that will
generalize to the definition of the equivariant quantum LR
coefficients.
\subsection{Quantum cohomology}
 The (small) quantum
cohomology of $X$ is a graded, commutative $\z[q]-$algebra with
unit, where the complex degree of $q$ is equal to $m$. It has a
$\z[q]-$basis $\{\sigma_\lambda\}$ indexed by partitions $\lambda$
included in the $p \times (m-p)$ rectangle. By definition
$\sigma_\lambda$ has complex degree equal to $|\lambda|$. The
multiplication, denoted by $\star$ , is defined by:
\[ \sigma_\lambda \star
\sigma_\mu = \sum_{d \geqslant 0} \sum_{\nu}
c_{\lambda,\mu}^{\nu,d}q^d \cdot \sigma_{\nu} \] The sum is over
$\nu$ such that $|\lambda| + |\mu| = |\nu| + md$. The coefficients
$c_{\lambda,\mu}^{\nu,d}$ are called \textit{quantum
Littlewood-Richardson coefficients}, and they are a special case
of the (3-point, genus 0) \textit{Gromov-Witten invariants}. They
are equal to the number of rational curves of degree $d$ passing
through Schubert varieties
$\Omega_\lambda(F_\bullet)$,$\Omega_\mu(G_\bullet)$ and
$\Omega_{\nu^\vee}(H_\bullet)$, for general flags
$F_\bullet$,$G_\bullet$ and $H_\bullet$. It is a deep result that
this gives an {\it associative} operation. The degree $0$ quantum
LR-coefficients are equal to the ordinary LR-coefficients
$c_{\lambda,\mu}^\nu$. In other words, the quantum cohomology is a
deformation of the ordinary cohomology of X.

Denote the quantum cohomology of $X$ by $QH^\star(X)$. Further
study of this algebra was done by
\cite{W,Be,BCF,Bu1,Y,BKT1,FW,Po}. Recall a particular case of the
``quantum Pieri rule" first proved in \cite{Be}:
 \[\sigma_\lambda \star \sigma_\Box =
\sum_{\mu \rightarrow \lambda} \sigma_{\mu} + q \sigma_{\lambda^-}
\](the last term is omitted if $\lambda^-$ does not exist).  This formula
will be generalized to the equivariant setting in the next
section.

We recall next the formal definition of the quantum
LR-coefficients. Let $\mbotreix$ be the Kontsevich moduli space of
degree $d$ stable maps from (arithmetic) genus 0 rational curves
to $X$ with $3$ marked points (\cite{KM,FP}). Represent the closed
points of this space by $(C,p_1,p_2,p_3;f)$. There are evaluation
maps \[ ev_i : \mbotreig \longrightarrow Gr(p,m)\] which send a
stable map $(C,p_1,p_2,p_3;f)$ to $f(p_i)$ and the forgetful map
\[\pi:\mbotreig \longrightarrow \overline{\mathcal{M}}_{0,3}\simeq
pt \]

\begin{prop}\label{qdef} The quantum LR-coefficient $c_{\lambda,\mu}^{\nu,d}$ is
equal to
\[\pi_\star(ev_1^\star(\sigma_\lambda) \cup ev_2^\star(\sigma_\mu) \cup
ev_3^\star(\sigma_{\nu^\vee}))
\]where $\sigma_\alpha$ is the Schubert class defined in \S 2.1. and
$\pi_\star$ is the Gysin morphism (\S 2.1).\end{prop}
\begin{proof} See e.g. Lemma 14 in \cite{FP}. \end{proof}

If $d=0$ one gets $\mbotreig \simeq Gr(p,m)$, so the definition
above becomes the geometric definition of the classical
LR-coefficients. We will see that this definition of quantum LR
coefficients generalizes to the definition of the equivariant
quantum LR coefficients.

\section{Equivariant quantum Schubert calculus}

\subsection{Equivariant quantum cohomology}
  The definition of the equivariant quantum
cohomology of a variety $X$ with a G-action was given in
\cite{GK}. Computations and properties of this object can be found
in \cite{GK,AS,Kim1,Kim2,Kim3}.
\\ \indent We restrict ourselves to the case when $X=Gr(p,m)$ and
$G = T \simeq(\cx^\star)^m$, as in the previous sections.
Moreover, we work with \textit{small} equivariant quantum
cohomology obtained by considering a certain restriction of the
(big) equivariant
quantum product.\\
\indent We state some of the properties of $T-$equivariant quantum
cohomology of $X$ as found in \cite{Kim1,Kim2}. Recall that
$\Lambda=H^\star_T(pt)=\z[T_1,...,T_m]$.
\begin{enumerate} \item The equivariant quantum cohomology of $X$
 is a graded associative, commutative $\Lambda[q]$-algebra with
 unit.
 \item
 It has an additive $\Lambda[q]-$basis $\{\sigma_\lambda\}$ indexed by
 partitions $\lambda$ in the $p \times (m-p)$
 rectangle. \item The ring multiplication, denoted $\circ$, is given by:
  \[ \sigma_\lambda
\circ \sigma_\mu = \sum_{d\geqslant 0}\sum_{\nu}
c_{\lambda,\mu}^{\nu,d}q^d \cdot \sigma_{\nu} \] where
$c_{\lambda,\mu}^{\nu,d}$ are the \textit{equivariant quantum
Littlewood-Richardson coefficients} (EQLR). \end{enumerate}

The EQLR coefficients are a generalization of both the
(non-equivariant, pure) quantum and the equivariant LR
coefficients and have the following properties:

 (i) $c_{\lambda,\mu}^{\nu,d}$ is a homogeneous
polynomial in $\Lambda$ of degree $|\lambda|+|\mu|-|\nu|-md$ .

(ii) When $d=0$, $c_{\lambda,\mu}^{\nu,d}$ is equal to the
equivariant LR-coefficient $c_{\lambda,\mu}^{\nu}$ .

 (iii) If $|\lambda|+|\mu|=|\nu|+md$ (i.e. if
$c_{\lambda,\mu}^{\nu,d}$ has polynomial degree 0),
$c_{\lambda,\mu}^{\nu,d}$ is the quantum LR-coefficient.

Note that (ii) and (iii) imply that the equivariant quantum
cohomology algebra is a graded deformation of both equivariant and
quantum cohomology of $X$.

Kim's definition of the EQLR coefficients, adapted to our context,
will be given in the next section. Property (1) is proved in Prop.
\ref{eqqcoh} of the next section, and properties (ii) and (iii)
are respectively Claim 1 and Claim 2 within the proof of this
proposition (found in the Appendix). Property (i) will hold by the
definition of the EQLR coefficients. Also by definition (see e.g.
\cite{Kim1} \S 4, (iv)) the equivariant quantum cohomology is
isomorphic, as a $\Lambda[q]-$module, with the free module
$H^\star_T(X) \otimes \z[q]$. Property (2) is equivalent to this
fact.

\subsection{Equivariant Littlewood-Richardson coefficients}The goal
of this section is to present the definition of the EQLR
coefficients. The main references are \cite{Kim1,Kim2}.

Recall that $\sigma_\lambda^T$ (resp.
$\widetilde{\sigma}_\lambda^T$) denote the equivariant cohomology
Schubert classe defined with respect to the standard (resp. the
opposite) flag. Recall also the following diagram from Section
2.3:
$$ \begin{CD}
 \mbotreix @>{ev_i}>> X \\ @V {\pi} VV
\\ \overline{\mathcal{M}}_{0,3} \simeq pt  \end{CD} $$
where $ev_i$ is the evaluation at the $i-$th point, $1 \leqslant i
\leqslant 3$. The $T-$action on $X$ induces a $T-$action on
$\mbotreix$ by: \[ t \cdot (C,p_1,p_2,p_3;f) :=
(C,p_1,p_2,p_3;\tilde{f}) \] where $\tilde{f}(x):=t \cdot f(x)$,
for $x$ in $C$ and $t$ in $T$. All the maps involved in the above
diagram are $T-$equivariant, therefore they determine a diagram$$
\begin{CD}
 \mbotreix_{T}:= ET \times_T \mbotreix @>{ev_i^T}>> X_{T} \\ @V {\pi^T} VV
\\ ET \times_T \overline{\mathcal{M}}_{0,3} \simeq BT  \end{CD} $$
\indent The EQLR coefficients generalize both their equivariant
and quantum versions. Given the the definitions of the latter ones
in Sections 2.2 and 2.3 , there is only one sensible choice: \[
c_{\lambda,\mu}^{\nu,d} = \pi_\star^T((ev_1^T)^\star
(\sigma_{\lambda}^T) \cup (ev_2^T)^\star (\sigma_{\mu}^T) \cup
(ev_3^T)^\star (\widetilde{\sigma}_{\nu^\vee}^T))
\]where $\pi_\star^T$ is the equivariant Gysin morphism.

Following Kim, we define the equivariant quantum cohomology. Let
$(A, \circ)$ be the graded $\Lambda[q]$-module having a
$\Lambda[q]$-basis $\{\sigma_\lambda\}$ indexed by partitions
$\lambda$ included in the $p \times (m-p)$ rectangle. The degrees
are the usual ones (see \S 2 above). Define a multiplication,
denoted $\circ$, among the basis elements of $A$ as follows: \[
\sigma_\lambda \circ \sigma_\mu = \sum_{d \geqslant 0}\sum_{\nu}
c_{\lambda,\mu}^{\nu,d}q^d \cdot \sigma_{\nu} \]

\begin{prop}[\cite{Kim2},\cite{Kim1}]\label{eqqcoh}
$(A, \circ)$ is a commutative, associative $\Lambda[q]$-algebra
with unit. There are canonical isomorphisms \begin{enumerate}
\item $A/\langle \Lambda^+ \cdot A \rangle \simeq QH^\star(X)$ as
$\z[q]$-algebras \item $A/\langle q \cdot A \rangle \simeq
H^\star_T(X)$ as $\Lambda$-algebras.
\end{enumerate} sending a basis element $\sigma_\lambda$ to the corresponding
$\sigma_\lambda$ in $QH^\star(X)$, respectively to
$\sigma_\lambda^T$ in $H^\star_T(X)$. $\Lambda^+$ denotes the
ideal of elements in $\Lambda$ of (strictly) positive degree.
\end{prop}

\begin{proof} The proof is given in the Appendix. \end{proof}

\noindent {\bf Notation:} The algebra $A$ from Prop. \ref{eqqcoh}
is the $T-$equivariant quantum cohomology algebra
of $X$ and it is denoted by $\eqqhx$.\\

\noindent \textit{Remark:} There is another description of
$\eqqhx$, involving a presentation with generators and relations.
This presentation was first computed in \cite{GK} for
Grassmannians and complete flag manifolds, then in \cite{AS,Kim1}
for partial flag manifolds. These presentations were used to
derive presentations for the corresponding non-equivariant quantum
cohomology algebras. Also, equivariant quantum cohomology has been
successfully used to study Mirror Symmetry phenomena (see
\cite{G}).

\section{Proof of the Main Lemma}
In this section we prove the Main Lemma and an additional
vanishing result for the EQLR coefficients. For that, we need to
introduce some results due to A. Buch (see \cite{Bu1,BKT1}). Let
$f:(C,p_1,p_2,p_3) \longrightarrow X$ be a stable map of degree
$d$, where the curve $C$ is isomorphic to a tree of
$\mathbb{P}^1$'s. The \textit{kernel of $f$}, denoted $ker(f)$, is
the largest subspace that is contained in all $f(x)$ for $x \in
C$. Similarly, define the \textit{span of $f$}, denoted $span(f)$,
to be the smallest subspace that contains all $f(x)$ for $x \in
C$. \\

\noindent \textit{Remark:} The above definitions of the kernel and
the span are a slight generalization of the definitions in
\cite{Bu1} since we allow the curve $C$ to be reducible.

The following two results can be found in \cite{Bu1} when $C
\simeq \mathbb{P}^1$.

\begin{prop}\label{Buch1} The kernel of $f$ has dimension at least $p-d$ and
the span of $f$ has dimension at most $p+d$.
\end{prop} \begin{proof}  The curve $C=
\bigcup_{i=1}^s C_i$ is a tree of rational curves $C_i \simeq
\mathbb{P}^1$, and $f$ restricted to each $C_i$ has some degree
$d_i$ such that $\sum_{i=1}^s d_i = d$. We use induction on the
number $s$ of components of $C$. If s=1, the assertion is Lemma 1
in \cite{Bu1}. Suppose $s > 1$. Assume that $C = C^{(1)} \bigcup
C^{(2)}$, where $C^{(1)}$ and $C^{(2)}$ are trees of rational
curves in $C$ intersecting in some point $x \in C$. Let $d^{(i)}$
be the degree of $f$ restricted to $C^{(i)}$. Let $W^{(i)}$,
$K^{(i)}$ be the span respectively the kernel of $f$ restricted to
$C^{(i)}$. The induction hypothesis implies that $\dim(W^{(i)})
\leqslant p+ d^{(i)}$, and $\dim(K^{(i)}) \geqslant p- d^{(i)}$.
But $C^{(1)}$ and $C^{(2)}$ intersect in a unique point $x$, and
$f(x)$ is a space of dimension $p$. Thus both $W^{(1)}$ and
$W^{(2)}$ contain the space $f(x)$. It follows that $\dim( span f)
= \dim\langle W^{(1)} + W^{(2)} \rangle \leqslant p + d^{(1)} +
d^{(2)} = p + d $.

For the kernels, note that both $K^{(1)}$ and $K^{(2)}$ are
contained in $f(x)$, with codimensions at most $d^{(1)}$
respectively $d^{(2)}$. Then the codimension of their intersection
is at most $d^{(1)}+d^{(2)} = d$, which shows that $\dim (ker f) =
\dim(K^{(1)} \cap K^{(2)}) \geq p-d$.
\end{proof}

Denote by $\overline{\lambda}(d)$ the partition obtained from
$\lambda$ by removing its first $d$ rows and by $\hat{\lambda}(d)$
the partition obtained from $\lambda$ by removing the leftmost $d$
columns. To be precise, if $\lambda=(\lambda_1,...,\lambda_p)$
then $\overline{\lambda}(d)=(\lambda_{d+1},...,\lambda_p)$ while
the $i-$th part $(\hat{\lambda}(d))_i$ of $\hat{\lambda}(d)$ is
equal to $\max(\lambda_i-d,0)$. Recall that $\sigma_{(0)}=1$.
\begin{prop}\label{Buch2} Let $f:(C,p_1,p_2,p_3) \longrightarrow X$
be a stable map of degree $d$, let $K$ be a $(p-d)$-dimensional
subspace of the kernel of $f$ and let $W$ be a $(p+d)$-dimensional
subspace containing the span of $f$. For any complete flag
$F_\bullet: 0 \subset F_1 \subset ... \subset F_m=\cx^m$, if the
image of $f$ intersects $\Omega_\lambda(F_\bullet)$ then $K$
belongs to the Schubert variety
$\Omega_{\overline{\lambda}(d)}(F_\bullet)$ in $Gr(p-d,m)$ and $W$
belongs to the Schubert variety
$\Omega_{\hat{\lambda}(d)}(F_\bullet)$ in $Gr(p+d,m)$.
\end{prop}
\begin{proof} The Proposition is Lemma 2 in \cite{Bu1} for
$C \simeq \mathbb{P}^1$, but the proof for general $C$ is the
same.
\end{proof}
\begin{lemma}\label{vanishing} Let $\lambda, \mu, \nu$ be three partitions in the $p \times (m-p)$
rectangle such that one of the following holds:
\begin{enumerate} \item $d<p$ and $\sigma_{\overline{\lambda}(d)}
\cdot \sigma_{\overline{\nu^\vee}(d)}=0$ in $H^\star(Gr(p-d,m))$.

\item $d<m-p$ and $\sigma_{\hat{\lambda}(d)} \cdot
\sigma_{\widehat{\nu^\vee}(d)}=0$ in $H^\star(Gr(p+d,m))$.
\end{enumerate} Then $c_{\lambda,\mu}^{\nu,d}=0$.
\end{lemma}

To prove the Lemma we need the following fact: \\

\noindent {\bf Fact 1.} Let $F:X' \longrightarrow Y$ be a
$T-$equivariant morphism of two algebraic varieties, with $Y$
smooth. Let $V$ be a $T-$invariant subvariety of $Y$ of
codimension $c$ and let $[V]_T \in H^{2c}_T(Y)$ be the equivariant
cohomology class of $V$. Then the equivariant cohomology pull-back
$F^\star_T([V]_T)$ is equal to $0$ if $F^{-1}(V)$ is empty. \\

This follows from the observation that $F^\star_T([V]_T)$ is
supported on $F^{-1}(V)_T$. Details can be found in \cite{Mi}.

\begin{proof}[Proof of Lemma \ref{vanishing}] The idea of proof to
show that if $c_{\lambda,\mu}^{\nu,d}$ is different from $0$ then
the intersection \[ ev_1^{-1}(\Omega_\lambda) \cap
ev_2^{-1}(\Omega_\mu) \cap
ev_3^{-1}(\widetilde{\Omega}_{\nu^\vee})\] in $\mbotreix$ is
nonempty. Then use Prop. \ref{Buch2} to get a contradiction.

Assume that $c_{\lambda,\mu}^{\nu,d}$ is not $0$. By the
definition of the equivariant Gysin morphism, the EQLR coefficient
is equal to
\[ c_{\lambda,\mu}^{\nu,d} =
\pi_{\star}^T((ev_1^T)^\star (\sigma_{\lambda}^T) \cup
(ev_2^T)^\star (\sigma_{\mu}^T) \cup (ev_3^T)^\star
(\widetilde{\sigma}_{\nu^\vee}^T))
\] Apply Fact 1 with $X'=\mbotreix$,
$Y=X \times X \times X$, $V = \Omega_{\lambda} \times \Omega_{\mu}
\times \widetilde{\Omega}_{\nu^\vee}^T$ and $F = (ev_1, ev_2,
ev_3)$, where $T$ acts diagonally on $Y$. Note that $V$ is
$T-$invariant and that the pull-back of its equivariant cohomology
class in $H^\star_T(Y)$ satisfies
\[ F^\star_T([V]_T) = (ev_1^T)^\star (\sigma_{\lambda}^T) \cup
(ev_2^T)^\star (\sigma_{\mu}^T) \cup (ev_3^T)^\star
(\widetilde{\sigma}_{\nu^\vee}^T) \] (this follows from the finite
dimensional approximation approach, discussed in \S 2.2 above).
Thus the inverse image $F^{-1}(V)$, which is equal to the
intersection
\[ ev_1^{-1} (\Omega_\lambda) \cap ev_2^{-1} (\Omega_\mu) \cap ev_3^{-1}
(\widetilde{\Omega}_\nu^\vee)\] must be nonempty in $\mbotreix$.
This amounts to the existence of a stable map \[ f:(C,p_1,p_2,p_3)
\longrightarrow X
\] whose image intersects $\Omega_\lambda, \Omega_\mu$
and $\widetilde{\Omega}_{\nu^\vee}$.

Suppose $p < d$. Choose $K_f$ to be a $(p-d)-$dimensional subspace
of the kernel of $f$ (such a $K_f$ exists by Proposition
\ref{Buch1}). Proposition \ref{Buch2} implies that $K_f$ belongs
to $\Omega_{\overline{\lambda}(d)} \cap \Omega_{\overline{\mu}(d)}
\cap \widetilde{\Omega}_{\overline{\nu^\vee}(d)}$. In particular
the intersection $\Omega_{\overline{\lambda}(d)} \cap
\widetilde{\Omega}_{\overline{\nu^\vee}(d)}$ is nonempty. But it
is a general fact that two Schubert varieties defined with respect
to opposite flags are in general position (see e.g. \cite{F1},
pag. 149) . It follows that the cohomology product
$\sigma_{\overline{\lambda}(d)} \cdot
\sigma_{\overline{\nu^\vee}(d)}$ must be nonzero in
$H^\star(Gr(p-d,m))$, contradicting the hypothesis (1).

The case when $d < m-p$ is treated in a similar fashion, using a
$(p+d)-$dimensional space $W_f$ including the span of $f$.
\end{proof} \indent The key result of this paper is a sufficiently
general condition that implies one of the hypothesis of Lemma
\ref{vanishing}, therefore giving a sufficient condition for the
vanishing of the EQLR coefficients. This condition is spelled out
in the Main Lemma. We divide its proof into two other lemmas,
which we prove first, corresponding to the hypotheses (1) and (2)
of Lemma \ref{vanishing}.
\begin{lemma}\label{lemma A} Let $\lambda,\mu,\nu$ be three partitions included in
the $p \times (m-p)$ rectangle and let $d$ be a positive integer.
Suppose that $d <p$ and that $|\lambda|+d^2> |\nu|+ md$. Then
$\sigma_{\overline{\lambda}(d)} \cdot
\sigma_{\overline{\nu^\vee}(d)}=0$ in $H^\star(Gr(p-d,m))$.
\end{lemma}
\begin{proof} To prove the Lemma, it is enough to verify the following
inequality:
\[|\overline{\lambda}(d)|+|\overline{\nu^\vee}(d)|>(p-d)(m-p+d) = \dim (
Gr(p-d,m))\] Let $\lambda =(\lambda_1,...,\lambda_p)$ and
$\nu^\vee=(\rho_1,...,\rho_p)$. Then $\overline{\lambda}(d) =
(\lambda_{d+1},...,\lambda_p)$ and
$\overline{\nu^\vee}(d)=(\rho_{d+1},...,\rho_p)$. The fact that
$|\lambda|+d^2> |\nu|+ md$ implies that $d^2+|\lambda|+|\nu^\vee|
> p(m-p)+md$. Then \[ \sum_{i=1}^d \lambda_i + \sum_{i=d+1}^p
\lambda_i + \sum_{j=1}^d \rho_j+ \sum_{j=d+1}^p \rho_j >
p(m-p)+md-d^2 \] hence
\[ \sum_{i=d+1}^p (\lambda_i + \rho_{i}) > p(m-p)+md-d^2
-\sum_{i=1}^d \lambda_i- \sum_{j=1}^d \rho_j
\]But \[\sum_{i=1}^d \lambda_i+ \sum_{j=1}^d \rho_j \leq 2d(m-p) \]
hence
\[\sum_{i=d+1}^p (\lambda_i + \rho_{i})
> p(m-p)+md-d^2 - 2d(m-p) = (p-d)(m-p+d)\]which finishes the proof of the Lemma.
\end{proof}
\begin{lemma}\label{lemma B}Let $\lambda,\mu,\nu$ be three partitions included in
$p \times (m-p)$ rectangle and let $d>0$ be an integer. Suppose
that $d < m-p$ and that $|\lambda|+d^2> |\nu|+ md$. Then
$\sigma_{\hat{\lambda}(d)} \cdot \sigma_{\widehat{\nu^\vee}(d)}=0$
in $H^\star(Gr(p+d,m))$.
\end{lemma}
\begin{proof} Consider the conjugation isomorphism $\psi : Gr(m-p,m)
\longrightarrow Gr(p,m)$. $\psi$ induces an isomorphism
$\psi^\star: H^\star (Gr(p,m) \longrightarrow H^\star(Gr(m-p,m))$
and it is well-known that $\psi^\star$ sends the class
$\sigma_\lambda \in H^\star(Gr(p,m))$ to the class
$\sigma_{\lambda'} \in H^\star(Gr(m-p,m))$, where $\lambda'$ is
the partition conjugate to $\lambda$.

Note that the hypotheses of the Lemma \ref{lemma A} are satisfied
when using partitions $\lambda', \mu', \nu'$ in the $(m-p) \times
p $ rectangle. Therefore $\sigma_{\overline{\lambda'}(d)} \cdot
\sigma_{\overline{(\nu')^\vee}(d)}=0$ in $H^\star(Gr(m-p-d,m))$.
Using again the conjugation isomorphism for the Grassmannians
$Gr(m-p-d,m)$ and $Gr(p+d,m)$ gives that
\[\sigma_{{\overline{\lambda'}(d)}'} \cdot
\sigma_{{\overline{(\nu')^\vee}(d)}'}=0\] in $H^\star(Gr(p+d,m))$.
Note that $(\nu')^\vee = (\nu^\vee)'$. To finish the proof it is
enough to prove the following combinatorial fact: \\
\textit{Fact :} Let $a,b,d$ be positive integers such that $d < b
$ and let $\lambda$ be a partition included in the $a \times b$
rectangle. Then $(\overline{\lambda'}(d))' = \hat{\lambda}(d)$ in
the $(a+d)\times(b-d)$ rectangle. \\ \textit{Proof:} Let $\lambda'
= (t_1,...,t_b)$. Then
$(\hat{\lambda}(d))'=(t_{d+1},...,t_b)=\overline{\lambda'}(d) $,
which implies that $\hat{\lambda}(d) = ((\hat{\lambda}(d))')'
=(\overline{\lambda'}(d))'$. \end{proof}
Concluding, the two previous lemmas add up to: \\

\noindent \textbf{Main Lemma.} \textit{ Let $\lambda,\mu,\nu$ be
three partitions included in $p \times (m-p)$ rectangle and let
$d>0$ be an integer. Suppose that $|\lambda|+d^2> |\nu|+ md$. Then
$c_{\lambda,\mu}^{\nu,d}=0$.}
\begin{proof}[Proof of the Main Lemma:]
\indent First note that $d \neq p$. Indeed, $|\lambda|+d^2> |\nu|+
md$ and $\lambda \subset (m-p)^p$ implies that $p(m-p)+d^2
> |\nu|+ md$. $p=d$ would imply $d^2 - d^2>|\nu|$ which is
impossible. It remains to study the cases when $d < p $ and $d>p$.

If $d < p $, apply Lemma \ref{lemma A} to get that
$\sigma_{\overline{\lambda}(d)} \cdot
\sigma_{\overline{\nu^\vee}(d)}=0$ in $H^\star(Gr(p-d,m))$. Then
Lemma \ref{vanishing}, statement (1), implies that
$c_{\lambda,\mu}^{\nu,d}=0$.

If $d > p$, we claim that $d < m-p$. Indeed, since
$|\nu^\vee|=p(m-p)-|\nu|$, we can rewrite the inequality
$|\lambda| + |\mu| \geqslant |\nu| + md$ as
\[|\lambda| + |\mu| + |\nu^\vee| \geqslant p(m-p) + md\] But $|\lambda|\leqslant p(m-p)
< d(m-p)$, hence
\[ d(m-p) + |\mu| + |\nu^\vee| > |\lambda|+ |\mu| +
|\nu^\vee| \geqslant p(m-p)+md.\] Then \[|\mu| + |\nu^\vee| >
p(m-p) +pd.
\] Since $|\mu|,|\nu^\vee| \leqslant p(m-p)$, one has that $pd < p(m-p)$,
which implies $d < m-p$, as needed. Hence we can apply Lemma
\ref{lemma B} to get that $\sigma_{\hat{\lambda}(d)} \cdot
\sigma_{\widehat{\nu^\vee}(d)} = 0$ in $H^\star(Gr(p+d,m))$.
Finally, Lemma \ref{vanishing}, statement (2) gives that
$c_{\lambda,\mu}^{\nu,d}=0$.
\end{proof}
\indent An immediate application of the Main Lemma is the next
Corollary, which shows the vanishing of the mixed EQLR Pieri
coefficients.

\begin{cor}\label{Pieri vanishing} Let $\lambda,\nu$ be two partitions included in
$p \times (m-p)$ rectangle and let $d$ be a positive integer. Then
$c_{\lambda,\Box}^{\nu,d}=0$, unless $d=1$ and $\nu =\lambda^-$.
\end{cor}
\begin{proof} If the polynomial degree of $c_{\lambda,\Box}^{\nu,d}$ is equal
to $0$, the assertion follows from the quantum Pieri rule (\S
2.3). If the polynomial degree of  $c_{\lambda,\Box}^{\nu,d}$ is
positive, since the partition $\Box$ is included in the $d \times
d$ square, note that
\[ |\lambda|+ d^2 \geqslant |\lambda|+|\Box| > |\nu| + md \]
so the conclusion follows from the Main Lemma.
\end{proof}

\noindent {\it Remark:} There is another proof of this result,
which doesn't use the Main Lemma, and which generalizes to any
homogeneous space $G/P$. However, this proof is weaker, in the
sense that it does not imply the vanishing result from the Main
Lemma. That is why we have chosen the proof above. \\

Next we prove another vanishing result, to be used later (\S 7).
\begin{prop}\label{full box vanishing} Let $\lambda$ be a partition included in the $p
\times (m-p)$ rectangle. Then $c_{\lambda,(m-p)^p}^{(0),d}=0$ for
$d < min\{p,m-p\}$. \end{prop} \begin{proof} By Lemma
\ref{vanishing}, the result follows if
$\sigma_{\overline{(m-p)^p}(d)} \cdot
\sigma_{\overline{(m-p)^p}(d)}=0$ in $H^\star(Gr(p-d,m))$. For
that, it is enough to note that
\[|\overline{(m-p)^p}(d)|+|\overline{(m-p)^p}(d)|=2(m-p)(p-d)>(p-d)(m-p+d)=\dim
(Gr(p-d,m))\] where the last inequality follows from the
assumption $d< m-p$.
\end{proof}
\section{Equivariant quantum Pieri rule} We prove the equivariant
quantum Pieri rule, then a recursive formula for the EQLR
coefficients. These will be the main ingredients of an algorithm
to compute the EQLR coefficients (cf. \S 7). Recall that the
equivariant quantum cohomology of $X$ is denoted by $\eqqhx$ and
that the equivariant coefficient $c_{\lambda,\Box}^\lambda$ is
given by
\[c_{\lambda,\Box}^\lambda =\sum_{i \in I(\lambda)} T_i - \sum_{j=m-p+1}^m T_j\]
where $I(\lambda)$ encodes the positions of the vertical steps in
the partition $\lambda$ (Prop. \ref{equivariant reccurence}).

\begin{thm}[equivariant quantum Pieri rule]\label{EQ Pieri} The following formula
holds in $\eqqh$:
 \begin{equation}\label{eq Pieri}\sigma_\lambda \circ \sigma_\Box =
\sum_{\mu \rightarrow \lambda} \sigma_{\mu} + c_{\lambda,\Box
}^\lambda \sigma_\lambda + q \sigma_{\lambda^-} \end{equation}
where the last term is omitted if $\lambda^-$ does not exist.
\end{thm}
\begin{proof} The equivariant quantum Pieri rule is a deformation
of both equivariant and quantum Pieri rules. In particular, it
must contain at least the terms on the right side. It remained to
prove that there are not any other terms. The only possibilities
are terms of the form $c_{\lambda,\Box}^{\nu,d} \sigma_\nu$ where
$d>0$ \textit{and} the polynomial degree of
$c_{\lambda,\Box}^{\nu,d}$ is positive. But Corollary \ref{Pieri
vanishing} shows that in this case $c_{\lambda,\Box}^{\nu,d}=0$,
as claimed.
\end{proof}

The equivariant quantum Pieri rule, commutativity and an
associativity relation implies a quantum generalization of an
equation satisfied by the equivariant LR coefficients (cf. Prop.
2.1 or \cite{MS,O,KT1}). It relates the EQLR coefficient
$c_{\lambda,\mu}^{\nu,d}$ with coefficients of degree $d$ one
smaller and coefficients of polynomial degree one larger.

\begin{prop}\label{reccurence} The EQLR coefficients satisfy the following
equation:
\begin{equation}\label{*} (\sum_{i \in
I(\nu)} T_i - \sum_{j \in I(\lambda)} T_j)\cdot
c_{\lambda,\mu}^{\nu,d} = \sum_{\delta \rightarrow \lambda}
c_{\delta,\mu}^{\nu,d} - \sum_{\nu \rightarrow \zeta}
c_{\lambda,\mu}^{\zeta,d} + c_{\lambda^- ,\mu}^{\nu,d-1} -
c_{\lambda,\mu}^{\nu^+,d-1} \end{equation} for any partitions
$\lambda,\mu,\nu$ and any nonnegative integer $d$. As usual, the
third (resp. the fourth) term in the right side is omitted if
$\lambda^-$ (resp. $\nu^+$) does not exist in the $p \times (m-p)$
rectangle. Both these terms are omitted if $d=0$.
\end{prop}

\begin{proof} We use the EQ Pieri rule and the associativity relation
$\sigma_\Box \circ (\sigma_\lambda \circ \sigma_\mu) =
(\sigma_\Box \circ \sigma_\lambda) \circ \sigma_\mu$. Recall that
$\sigma_\lambda \circ \sigma_\Box = \sigma_\Box \circ
\sigma_\lambda$ for any partition $\lambda$. The result will
follow from the identification of the coefficient of $q^d
\sigma_\nu$ in both sides of the relation. Indeed, we have
\begin{eqnarray*}\sigma_\Box \circ (\sigma_\lambda \circ \sigma_\mu)
 &=& \sigma_\Box \circ ( \sum_{d, \rho} q^d c_{\lambda,\mu}^{\rho,
d} \sigma_{\rho}) \\ &=& \sum_{d, \rho} q^d
c_{\lambda,\mu}^{\rho,d} (\sum_{\theta \rightarrow \rho}
\sigma_{\theta} + q \sigma_{\rho^- }+ c_{\Box,\rho}^{\rho}
\sigma_{\rho})\end{eqnarray*} and
\begin{eqnarray*}(\sigma_\Box \circ \sigma_\lambda) \circ
\sigma_\mu &=& (\sum_{\delta \rightarrow \lambda} \sigma_{\delta}
+ q \sigma_{\lambda^-} + c_{\Box,\lambda}^{\lambda}\sigma_\lambda)
\circ \sigma_{\mu} \\ & = & \sum_{\delta \rightarrow \lambda}
(\sum_{\alpha,d_1} q^{d_1}c_{\delta,\mu}^{\alpha,d_1}
\sigma_\alpha) + q (\sum_{\beta,d_2}
q^{d_2} c_{\lambda^-,\mu}^{\beta, d_2} \sigma_{\beta}) \\
& & + c_{\Box,\lambda}^{\lambda}(\sum_{d_3, \gamma} q^{d_3}
c_{\lambda,\mu}^{\gamma,d_3}\sigma_{\gamma})
\end{eqnarray*}
Collecting the coefficient of $q^d \sigma_\nu$ from both sides
gives
\begin{equation}\label{local} \sum_{\nu \rightarrow \zeta} c_{\lambda,\mu}^{\zeta,d}+
c_{\lambda,\mu}^{\nu^+,d-1} + c_{\Box,\nu}^\nu
c_{\lambda,\mu}^{\nu,d} = \sum_{\delta \rightarrow \lambda}
c_{\delta,\mu}^{\nu,d} + c_{\lambda^-,\mu}^{\nu,d-1} +
c_{\Box,\lambda}^{\lambda}c_{\lambda,\mu}^{\nu,d} \end{equation}
Note that the difference $c_{ \Box, \nu }^\nu - c_{ \Box, \lambda
}^\lambda$ is equal to
\[  c_{ \Box, \nu }^\nu - c_{ \Box, \lambda }^\lambda = \sum_{i \in
I(\nu)} T_i - \sum_{j \in I(\lambda)} T_j. \] Then formula
(\ref{*}) follows by rearranging the terms of expression
(\ref{local}).
\end{proof}

\section{Two formulae}
In this section we prove two formulae for some special classes of
EQLR coefficients, used in the algorithm computing these
coefficients (see \S 7). All the results from now on will be
algorithmic, and the EQLR coefficients are no longer considered to
be homogeneous polynomials, but (possibly nonhomogeneous) {\em
rational functions} in the fraction field of $\Lambda$. Thus, a
priori, we may have a nonzero coefficient
$c_{\lambda,\mu}^{\nu,d}$ with negative
$|\lambda|+|\mu|-|\nu|-md$. The latter {\em quantity} will still
be referred to as ``polynomial degree'', to distinguish it from
the ``degree'' $d$ of $c_{\lambda,\mu}^{\nu,d}$. Unless otherwise
stated, the only assumption about these coefficients is that they
satisfy the following equation, obtained by rewriting the equation
(\ref{*}) above as:
\begin{equation}\label{**}c_{\lambda,\mu}^{\nu,d} =
\frac{\sum_{\delta \rightarrow \lambda}
c_{\delta,\mu}^{\nu,d}}{F_{\nu,\lambda}} - \frac{\sum_{\nu
\rightarrow \zeta} c_{\lambda,\mu}^{\zeta,d}}{F_{\nu,\lambda}} +
\frac{c_{\lambda^- ,\mu}^{\nu,d-1} -
c_{\lambda,\mu}^{\nu^+,d-1}}{F_{\nu,\lambda}}
\end{equation} for $\nu$ different from $\lambda$.
$F_{\nu,\lambda}$ denotes the polynomial
\[F_{\nu,\lambda} = \sum_{i \in I(\nu)} T_i - \sum_{j \in I(\lambda)} T_j. \]
Recall that we omit the coefficient $c_{\lambda^- ,\mu}^{\nu,d-1}$
(resp. $c_{\lambda,\mu}^{\nu^+,d-1}$) in the third term of the
right side of (\ref{**}) if $\lambda^-$ (resp. $\nu^+$) does not
exist. The entire third term is omitted if $d=0$. The propositions
below show that the EQLR coefficient $c_{\lambda,\mu}^{\nu,d}$
(for some special $\lambda,\mu,\nu$) is determined by EQLR
coefficients of degree $d-1$, and possibly some coefficients of
degree $d$. Their proof is by induction on $|\lambda|-|\nu|$.

\begin{prop}\label{inclusion} Let $\lambda,\mu,\nu$ be partitions such that $\lambda$ is
not included in $\nu$ and let $d$ be a nonnegative integer. Then
\begin{equation}\label{einclusion} c_{\lambda,\mu}^{\nu,d}= E_{\lambda,\mu,\nu}(d) \end{equation} where
$E_{\lambda,\mu,\nu}(d)$ is a linear homogeneous expression in
EQLR coefficients of degree $d-1$ with coefficients in
$R(\Lambda)$, the fraction field of $\Lambda(=\z[T_1,...,T_m])$.
If $d=0$ then $c_{\lambda,\mu}^{\nu,0}=E_{\lambda,\mu,\nu}(0)=0$.
\end{prop}

For the next proposition, let $\alpha$ and $\lambda$ be two
partitions such that $\alpha$ is included in $\lambda$. Define a
rational function $R_{\lambda,\alpha}$ in $R(\Lambda)$ as follows:

$$ R_{\lambda,\alpha}= \left \{
\begin{array}{ll}
 \sum \prod_{i=0}^{l-1}
\frac{1}{F_{\lambda,\alpha^{(i)}}} \ & \textrm{if } \lambda \neq \alpha\\
1 & \textrm {if } \alpha = \lambda
\end{array} \right.  $$

\noindent In the case $\lambda \neq \alpha$, $l$ denotes the
nonnegative integer $|\lambda|-|\alpha|$, and the sum  is over all
chains of partitions
\[ \lambda=\alpha^{(l)} \to \alpha^{(l-1)} \to ...\to \alpha^{(1)} \to
\alpha^{(0)}=\alpha. \]

\begin{prop}\label{equality} The EQLR coefficient $c_{\alpha,\lambda}^{\lambda,d}$
satisfies the following formula:
\begin{equation}\label{equation} c_{\alpha,\lambda}^{\lambda,d} =
R_{\lambda,\alpha}c_{\lambda,\lambda}^{\lambda,d} +
E'_{\lambda,\alpha}(d) \end{equation} where
$E'_{\lambda,\alpha}(d)$ is an $R(\Lambda)-$linear homogeneous
expression in EQLR coefficients of degree $d-1$. If $d=0$ then
$E'_{\lambda,\alpha}(0)=0$. Moreover, for any such $\lambda$ and
$\alpha$, $R_{\lambda,\alpha}$ is different from $0$.
\end{prop}

The proof of these propositions requires some notations. Let
$(\lambda,\nu)$ and $(\delta,\zeta)$ be two pairs of partitions
included in the $p \times (m-p)$ rectangle. Define
\[ (\lambda,\nu) <_1 (\delta,\zeta)\] if  $\delta \rightarrow
\lambda$ and $\nu = \zeta$ and
\[ (\lambda,\nu) <_2 (\delta,\zeta)\] if $\delta =
\lambda$ and $\nu \rightarrow \zeta$. With these notations, the
first two terms of the right side of equation (\ref{**}), when
this is applied to $c_{\lambda,\mu}^{\nu,d}$, contain coefficients
of the form $c_{\delta,\mu}^{\zeta,d}$ with either $(\lambda,\nu)
<_1 (\delta,\zeta)$ or $(\lambda,\nu) <_2 (\delta,\zeta)$. Note
also that if $(\lambda,\nu) <_i (\delta,\zeta)$ ($i=1,2$) then the
difference $|\delta|-|\zeta|$ is one larger than the difference
$|\lambda|-|\nu|$.

\begin{proof}[Proof of Proposition \ref{inclusion}] Use descending induction on the
difference $|\lambda|-|\nu|\leqslant p(m-p)$. If $|\lambda|-|\nu|=
p(m-p)$, then $\lambda = (m-p)^p$ (the partition having $p$ parts
of length $m-p$) and $\nu=(0)$. The first two terms of the
equation (\ref{**}) vanish when it is applied to
$c_{\lambda,\mu}^{\nu,d}=c_{(m-p)^p,\mu}^{(0),d}$, so
$c_{\lambda,\mu}^{\nu,d}$ is equal to $E_{(m-p)^p,\mu,(0)}(d)$
where
\[ E_{(m-p)^p,\mu,(0)}(d) = \frac{1}{F_{(0),(m-p)^p}}
(c_{(m-p-1)^{p-1},\mu}^{(0),d-1}-
c_{(m-p)^{p},\mu}^{(m-p,1^{p-1}),d-1}).
\] Here $(m-p,1^{p-1})$ is the partition having the first part equal to
$m-p$ and the next $p-1$ parts equal to $1$. The base of the
induction is proved.

Let $|\lambda|-|\nu|< p(m-p)$ such that $\lambda$ is not included
in $\nu$. Applying (\ref{**}) to $c_{\lambda,\mu}^{\nu,d}$ yields
coefficients of degree $d-1$ and coefficients
$c_{\delta,\mu}^{\zeta,d}$ with $(\lambda,\nu) <_i (\delta,\zeta)$
($i=1,2$). It is enough to show that each such coefficient
$c_{\delta,\mu}^{\zeta,d}$ is equal to a $R(\Lambda)-$linear
homogeneous expression $E_{\delta,\mu,\zeta}(d)$ in coefficients
of degree $d-1$. Note that $\lambda \subset \delta$ and $\zeta
\subset \nu$, and, since $\lambda$ is not included in $\nu$, it
follows that $\delta$ is not included in $\zeta$. Since
$|\delta|-|\zeta|=|\lambda|-|\nu|+1$, the induction hypothesis,
applied to $c_{\delta,\mu}^{\zeta,d}$, concludes the proof if
$d>0$. The same proof works for $d=0$, yielding now
$c_{\lambda,\mu}^{\nu,0}=0$, since in this case the last term of
(\ref{**}) is ignored. \end{proof}

\begin{proof}[Proof of Prop. \ref{equality}] Note that $\alpha$ is included in $\lambda$
implies that $|\lambda|-|\alpha| \geqslant 0$. We induct on
$|\lambda|-|\alpha|$. If $|\lambda|-|\alpha|=0$ then
$\alpha=\lambda$, and there is nothing to prove (in this case
$E'_{\lambda,\alpha}(d)=0$). Let $|\lambda|-|\alpha|>0$. In
particular $\lambda$ is not equal to $\alpha$, so one can apply
equation (\ref{**}) to $c_{\alpha,\lambda}^{\lambda,d}$. The first
two terms of the RHS of (\ref{**}) contain coefficients
$c_{\delta,\lambda}^{\zeta,d}$ with $(\alpha,\lambda)<_1
(\delta,\zeta)$ (first term) respectively $(\alpha,\lambda)<_2
(\delta,\zeta)$ (second term).

Consider first a coefficient $c_{\delta,\lambda}^{\zeta,d}$ from
the second term. Then $\delta=\alpha$ and $\lambda \rightarrow
\zeta$. In particular $\lambda$ is not included in $\zeta$, thus,
by formula (\ref{einclusion})
\begin{equation}\label{4}
c_{\alpha,\lambda}^{\zeta,d}=E_{\alpha,\lambda,\zeta}(d).
\end{equation} Consider now a coefficient $c_{\delta,\lambda}^{\zeta,d}$
from the first term. Then $\lambda=\zeta$ and $\delta \rightarrow
\alpha$. If $\delta$ is not included in $\lambda$, then
\begin{equation}\label{6}
c_{\delta,\lambda}^{\zeta,d}=c_{\delta,\lambda}^{\lambda,d}=E_{\delta,\lambda,\lambda}(d)
\end{equation} again by formula (\ref{einclusion}). If $\delta$ is included in
$\lambda$, by induction hypothesis
\begin{equation}\label{5}
c_{\delta,\lambda}^{\zeta,d}=c_{\delta,\lambda}^{\lambda,d}=
R_{\lambda,\delta} c_{\lambda,\lambda}^{\lambda,d} +
E'_{\lambda,\delta}(d) \end{equation} Combining
(\ref{**}),(\ref{4}),(\ref{6}),(\ref{5}) and noting that
\[ R_{\lambda,\alpha} = \frac{1}{F_{\lambda,\alpha}}
\sum_{} R_{\lambda,\delta} \] (the sum is over all $\delta$ such
that $\delta \to \alpha$ and $\delta$ included in $\lambda$)
yields
\begin{eqnarray*} c_{\alpha,\lambda}^{\lambda,d}& =&
\frac{1}{F_{\lambda,\alpha}}\Bigl(\sum_{} \bigl(R_{\lambda,\delta}
c_{\lambda,\lambda}^{\lambda,d} + E'_{\lambda,\delta}(d)\bigr) +
\sum E_{\delta,\lambda,\lambda}(d)\Bigr) +\\& &
\frac{1}{F_{\lambda,\alpha}}\Bigl(\sum_{\lambda \rightarrow
\zeta}E_{\alpha,\lambda,\zeta}(d)\Bigr)+
 \frac{c_{\lambda^- ,\mu}^{\nu,d-1}
- c_{\lambda,\mu}^{\nu^+,d-1}}{F_{\nu,\lambda}} = \\ &= &
R_{\lambda,\alpha} c_{\lambda,\lambda}^{\lambda,d} +
E'_{\lambda,\alpha}(d). \end{eqnarray*} The first (resp. the
second) sum is over all $\delta$ such that $\delta \to \alpha$ and
included (resp. not included) in $\lambda$, while
$E'_{\lambda,\alpha}(d)$ is obtained by collecting all the terms
involving coefficients of degree $d-1$. If $d=0$ the same proof
shows that $E'_{\lambda,\alpha}(0)=0$.

To finish the proof it remained to show that $R_{\lambda,\alpha}$
is not equal to zero. Since the partition $\alpha^{(j)}$ is
included in $\lambda$ ($j=0,...,|\lambda|-|\alpha|$) ,
$F_{\lambda,\alpha^{(j)}}$ is a linear homogeneous polynomial in
variables $T_1-T_2,...,T_{m-1}-T_m$ with {\it positive}
coefficients. This shows that there cannot be cancellations in the
sum defining $R_{\lambda,\alpha}(T)$.
\end{proof}

\section{An algorithm to compute the EQLR coefficients} Granting
the fact that the equivariant quantum multiplication is a
commutative operation with unit, the next theorem shows, by an
algorithm, that all EQLR coefficients $c_{\lambda,\mu}^{\nu,d}$
can be recovered from the Pieri ones and the recursive formula
(\ref{**}). The proof is by double induction, on the degree $d$,
then on the polynomial degree. The idea is to show, by descending
induction on polynomial degree, that equation (\ref{**}) implies
that $c_{\lambda,\mu}^{\nu,d}$ is determined by coefficients of
the form $c_{\alpha,\beta}^{\gamma,d-1}$, known by induction on
$d$, and by coefficients $c_{\chi,\chi}^{\chi,d}$ which are known
again by induction on $d$ (see formula (\ref{equation})).
\begin{thm}\label{recalg} The EQLR coefficients are determined (algorithmically)
by the following formulae:

\noindent (i) $c_{(0),(0)}^{(0),d}=0$ unless $d=0$, when it is
equal to $1$.

\noindent (ii) (commutativity)
$c_{\lambda,\mu}^{\nu,d}=c_{\mu,\lambda}^{\nu,d}$ for all
partitions $\lambda,\mu$ and $\nu$

\noindent (iii)(EQ Pieri) The coefficients
$c_{\Box,\lambda}^{\nu,d}$ from equation (\ref{eq Pieri}), for all
partitions $\lambda,\nu$.

\noindent (iv) Equation (\ref{**}), for all partitions
$\lambda,\mu,\nu$ such that $\lambda$ is different from $\nu$.
\end{thm}

Before proving the theorem, we would like to emphasize the
following corollary, which is a weaker, but useful, version of the
opening paragraph of this section:
\begin{cor} Let $(A, \diamond)$ be a graded, commutative,
associative $\Lambda[q]-$algebra with unit such that:

\noindent 1. $A$ has an additive $\Lambda[q]-$basis
$\{t_\lambda\}$ (graded as usual).

\noindent 2. The equivariant quantum Pieri holds, i.e.
\[t_\lambda \diamond t_\Box =
\sum_{\mu \rightarrow \lambda} t_{\mu} + c_{\lambda,\Box }^\lambda
t_\lambda + q t_{\lambda^-} \] where the last term is omitted if
$\lambda^-$ does not exist.

Then $A$ is canonically isomorphic to $\eqqh$, as
$\Lambda[q]-$algebras.
\end{cor}

\begin{proof}[Proof of the Corollary] The structure constants of
$A$ clearly satisfy (i)-(iii); (iv) follows from the associativity
of $A$ and the EQ Pieri rule (cf. Prop. \ref{reccurence}).
\end{proof}

\begin{proof}[Proof of Theorem \ref{recalg}] The algorithm has
three steps, with the main induction on $d$, the degree of the
EQLR coefficient $c_{\lambda,\mu}^{\nu,d}$. The base case ($d=0$)
is treated in Step 1:

\noindent Step 1: Compute the coefficient
$c_{\lambda,\mu}^{\nu,0}$ for all partitions $\lambda,\mu,\nu$.

In fact we prove that the hypotheses of the theorem imply that
$c_{\lambda,\mu}^{\nu,0}$ is equal to the equivariant coefficient
$c_{\lambda,\mu}^\nu$. The latter can then be computed using Prop.
\ref{equivariant reccurence}, by induction on $|\nu|-|\lambda|$
(see also \cite{KT1} \S 3, Cor. 1). To prove the equality of the
two coefficients it is enough to show that the coefficients
$c_{\lambda,\mu}^{\nu,0}$ satisfy the following formulae (see
Prop. \ref{equivariant reccurence}):

(a') $ c_{\lambda,\Box }^{\lambda,0} = \sum_{i \in I(\lambda)} T_i
- \sum_{j=m-p+1}^m T_j $

(b') $c_{\lambda,\lambda}^{\lambda,0} =\prod_{i \in I(\lambda), j
\in J(\lambda), i<j} (T_i - T_j)$, where $J(\lambda)$ encodes the
position of the horizontal steps in $\lambda$.

(c') $(c_{\lambda,\Box }^{\lambda} - c_{\mu,\Box }^\mu)\cdot
c_{\lambda,\mu}^{\lambda,0} = \sum_{\delta \rightarrow \mu}
c_{\lambda,\delta}^{\lambda,0}$ for any $\lambda, \mu$ such that
$\lambda \neq \mu$.

(d') $(c_{\nu,\Box }^\nu - c_{\lambda,\Box }^\lambda)\cdot
c_{\lambda,\mu}^{\nu,0} = \sum_{\delta \rightarrow \lambda}
c_{\delta,\mu}^{\nu,0} - \sum_{\nu \rightarrow \zeta}
c_{\lambda,\mu}^{\zeta,0}$ for any $\lambda,\mu,\nu$ such that
$\lambda \neq \nu$.

\begin{proof} Formula (a') follows from the equivariant quantum
Pieri (EQ Pieri) rule (iii); equation (d') follows from (iv) since
the last term of equation (\ref{**}) is omitted in this case.
Equation (c') follows from (d') (using commutativity (ii)) once we
show that $c_{\lambda,\mu}^{\nu,0}$ vanishes if $\lambda$ is not
included in $\nu$. This holds by Prop. \ref{inclusion}.

It remains to prove formula (b'). If $\lambda$ is equal to $(0)$
or $\Box$, this follows respectively from (i) and (ii). For a
bigger $\lambda$, equation (\ref{equation}) in Prop.
\ref{equality} implies that
\[ c_{\lambda,\lambda}^{\lambda,0} =
\frac{1}{R_{\lambda,\alpha}}c_{\alpha,\lambda}^{\lambda,0} \] for
any partition $\alpha$ included in $\lambda$. Note that the same
equation holds for the equivariant LR coefficients
$c_{\lambda,\mu}^\nu$ (use induction on the difference
$|\lambda|-|\alpha| \geqslant 0$ and \ equation (c) in Prop. 2.1 -
same proof as for Prop. \ref{equality} above). Then
\[c_{\lambda,\lambda}^{\lambda,0} =
\frac{1}{R_{\lambda,\Box}}c_{\lambda,\Box}^{\lambda,0} =
\frac{1}{R_{\lambda,\Box}}c_{\lambda,\Box}^{\lambda}=
c_{\lambda,\lambda}^\lambda \] which concludes the proof.
\end{proof}

Let now $d>0$, and assume, by induction on $d$, that all EQLR
coefficients of degree $d-1$ are known. \\

\noindent Step 2: Compute $c_{\lambda,\lambda}^{\lambda,d}$ for
each partition $\lambda$.

If $\lambda$ is $(0)$ or $\Box$ this is given respectively in (i)
and (ii). Then let $\lambda$ be of weight at least 2. The EQLR
coefficient
$c_{\Box,\lambda}^{\lambda,d}=c_{\lambda,\Box}^{\lambda,d}$
vanishes by (iii). Then equation (\ref{equation}) implies that
\[c_{\lambda,\lambda}^{\lambda,d}=-\frac{E'_{\lambda,\Box}(d)
}{R_{\lambda,\Box}}\]
and $E'_{\lambda,\Box}(d)$ is known by induction on $d$. \\

\noindent Step 3: Compute all coefficients
$c_{\lambda,\mu}^{\nu,d}$.

Within the main induction on $d$, we use descending induction on
$|\lambda|+ |\mu|-|\nu|-md$, the polynomial degree of
$c_{\lambda,\mu}^{\nu,d}$. This degree is at most $2p(m-p)-md$, in
which case $\lambda$ and $\mu$ are both equal to $(m-p)^p$ and
$\nu = (0)$. The coefficient $c_{(m-p)^p,(m-p)^p}^{(0),d}$ is
known by induction on $d$ (by Prop. \ref{inclusion}, since
$\lambda$ is not included in $\nu$). The base of the induction is
proved.

Assume now that the polynomial degree of $c_{\lambda,\mu}^{\nu,d}$
is less than $2p(m-p)-md$. If $\lambda=\mu=\nu$ apply Step 2. If
not, since $c_{\lambda,\mu}^{\nu,d}$ is equal to
$c_{\mu,\lambda}^{\nu,d}$, we can assume that $\lambda$ is
different from $\nu$. Using (iv), write $c_{\lambda,\mu}^{\nu,d}$
as a combination of EQLR coefficients of polynomial degree one
larger (first and second term in the right side) and coefficients
of degree $d-1$ (third term), the latter ones known by induction
on $d$. The polynomial degree induction, applied to the
coefficients appearing in the first two terms of the right side of
(\ref{**}) finishes the algorithm. \end{proof}

For a more efficient algorithm one may include the following
vanishing properties of the EQLR coefficients:

(i) $c_{\lambda,\mu}^{\nu,d}=0$ if it has negative polynomial
degree, i.e. if $|\lambda|+|\mu|-|\nu|-md<0$.

(ii) $c_{\lambda,\mu}^{\nu,d}=0$ if $|\lambda|+ d^2 > |\nu|+ md$
or $|\mu|+ d^2 > |\nu|+ md$ (Main Lemma).

(iii)$c_{\lambda,(m-p)^p}^{(0),d}=0$ if $d < min\{p,m-p\}$ (Proposition
4.7).\\

There is also a method to reduce the number of computations needed in
Step 2 of the algorithm. This method can be applied to compute the
coefficients $c_{\lambda,\lambda}^{\lambda,d}$ with nonnegative
polynomial degree (equivalently, with $|\lambda| \geqslant md$), provided
one knows certain coefficients of polynomial degree $0$ (hence, in
geometric interpretation, some (pure) quantum coefficients). Indeed,
given a partition $\lambda$ of weight at least $md$, choose a partition
$\alpha=\alpha(\lambda)$ included in $\lambda$, of weight $|\alpha|=md$.
Note that $c_{\alpha,\lambda}^{\lambda,d}$ has polynomial degree $0$.
Then one can solve for $c_{\lambda,\lambda}^{\lambda,d}$ in the equation
(\ref{equation}) from Prop. \ref{equality}, using now
$c_{\alpha,\lambda}^{\lambda,d}$ on the left
side.\begin{footnote}[4]{Theoretical algorithms computing quantum LR
coefficients have been obtained in \cite{BKT1} (it uses the intersection
theory on 2-step flag manifolds), or \cite{BCF} (the rim-hook algorithm).
For explicit calculations we recommend Buch's Littlewood-Richardson
calculator, implemented in C, which can be found at
\texttt{http://home.imf.au.dk/abuch/lrcalc/}.}
\end{footnote} The author has implemented this algorithm in Maple 7.

Using again double induction on the degree $d$ and on the
polynomial degree implies that a coefficient
$c_{\lambda,\mu}^{\nu,d}$ obtained by the algorithm is in fact a
{\em homogeneous} rational function {\em in the variables
$T_1-T_2,...,T_{m-1}-T_m$} of degree $|\lambda|+|\mu|-|\nu|-md$.
This follows from the fact that the function $R_{\lambda,\alpha}$,
for any pair $\alpha \subset \lambda$ (cf. Prop. \ref{equality}),
and the coefficients in the hypothesis of the algorithm satisfy
this property, and that all the terms in equation (\ref{**}) have
the same polynomial degree. Unfortunately, it is not apparent from
the algorithm that $c_{\lambda,\mu}^{\nu,d}$ is a polynomial.
Nevertheless, this can be deduced from its ``geometric
realization''. We have just proved:
\begin{cor} The EQLR coefficients $c_{\lambda,\mu}^{\nu,d}$ are
homogeneous polynomials in variables
$T_1-T_2,...,T_{m-1}-T_m$.\begin{footnote}[5]{A geometric proof of this
result can be obtained by using the canonical action on $Gr(p,m)$ of the
maximal torus $T(m) \simeq (\cx^\star)^{m}/\cx^\star$ in $PGL(m)$.
Details can be found in \cite{Mi,Mi2}.}\end{footnote}
\end{cor}

\section{Examples/Final remarks} In this section we give
several values of the EQLR coefficients
$c_{\lambda,\lambda}^{\lambda,d}$ for Grassmannians $Gr(p,m)$,
with $p,m$ small, then the multiplication table for
$QH^\star_T(Gr(2,4))$. The section ends with some remarks about
related work in progress.

\subsection{The coefficients $c_{\lambda,\lambda}^{\lambda,d}$
for small Grassmannians} \indent

\begin{tabular}{|c|c|c|c|c|}
\hline
p & m & $\lambda$ & d & $c_{\lambda,\lambda}^{\lambda,d}$  \\
\hline
  2 & 4 & (2,2) & 1 &  0 \\ \hline
  2 &  5 & (3,2) & 1 & 1\\ \hline
  2 &  5 & (3,3) & 1 & 0 \\ \hline
  2 & 6 & (3,3) & 1 &  1 \\ \hline
  2 & 6 & (4,2) & 1 & 1\\ \hline
  2 & 6 & (4,3) & 1 & $T_1+T_2+T_3-T_4-T_5-T_6$\\ \hline
  2 & 6 & (4,4) & 1 &  0 \\ \hline
\end{tabular}
\begin{tabular}{|c|c|c|c|c|}
\hline
p & m & $\lambda$ & d & $c_{\lambda,\lambda}^{\lambda,d}$  \\
\hline
  3 & 6 & (3,3) & 1 & 0 \\ \hline
  3 & 6 & (3,2,1) & 1 & 2 \\ \hline
  3 & 6 & (2,2,2) & 1 & 0 \\ \hline
  3 & 6 & (3,3,1) & 1 & $T_1+T_2-T_3-T_4$\\ \hline
  3 & 6 & (3,2,2) & 1 & $T_3+T_4-T_5-T_6$ \\ \hline
  3 & 6 & (3,3,2) & 1 & $(T_1+T_2-T_5-T_6)^2$ \\ \hline
  3 & 6 & (3,3,3) & 1 & 0 \\ \hline
\end{tabular}
\subsection{Multiplication table for $QH^\star_T(Gr(2,4))$}

 \indent \\
 \indent $\sigma_{(1)} \circ \sigma_{(1)} = \sigma_{(2)}
 + \sigma_{(1,1)} + (T_2 - T_3) \sigma_{(1)} $ \\
\indent $\sigma_{(1)} \circ \sigma_{(2)} = \sigma_{(2,1)}
 + (T_1 - T_3) \sigma_{(2)} $ \\
\indent $\sigma_{(1)} \circ \sigma_{(1,1)} = \sigma_{(2,1)}
 + (T_2 - T_4) \sigma_{(1,1)} $ \\
\indent $\sigma_{(1)} \circ \sigma_{(2,1)} = \sigma_{(2,2)}
 + (T_1 - T_4) \sigma_{(2,1)}+ q $ \\
\indent $\sigma_{(1)} \circ \sigma_{(2,2)} =
(T_1+T_2-T_3-T_4)\sigma_{(2,2)} + q\sigma_{(1)} $ \\
\indent $\sigma_{(2)} \circ \sigma_{(2)} =
\sigma_{(2,2)} + (T_1-T_2)\sigma_{(2,1)}+ (T_1-T_2)(T_1-T_3)\sigma_{(2)} $ \\
\indent $\sigma_{(2)} \circ \sigma_{(1,1)} =
(T_1-T_4)\sigma_{(2,1)}+ q $ \\
\indent $\sigma_{(2)} \circ \sigma_{(2,1)} =
(T_1-T_4)\sigma_{(2,2)} + (T_1-T_2)(T_1-T_4)\sigma_{(2,1)}+
q\sigma_{(1)} +(T_1-T_2)q$ \\
\indent $\sigma_{(2)} \circ \sigma_{(2,2)} =
(T_1-T_4)(T_1-T_3)\sigma_{(2,2)} + q\sigma_{(1,1)}
+(T_1-T_3)q\sigma_{(1)}$ \\
\indent $\sigma_{(1,1)} \circ \sigma_{(1,1)} =
\sigma_{(2,2)} + (T_3-T_4)\sigma_{(2,1)}+ (T_2-T_4)(T_3-T_4)\sigma_{(1,1)} $ \\
\indent $\sigma_{(1,1)} \circ \sigma_{(2,1)} =
(T_1-T_4)\sigma_{(2,2)} + (T_1-T_4)(T_3-T_4)\sigma_{(2,1)}+
q\sigma_{(1)} +(T_3-T_4)q$ \\
\indent $\sigma_{(1,1)} \circ \sigma_{(2,2)} =
(T_1-T_4)(T_2-T_4)\sigma_{(2,2)} + q\sigma_{(2)}
+(T_2-T_4)q\sigma_{(1)}$

$\sigma_{(2,1)} \circ \sigma_{(2,1)} =  (T_1-T_4)^2\sigma_{(2,2)}
+ (T_1-T_2)(T_1-T_4)(T_3-T_4)\sigma_{(2,1)} + q \sigma_{(2)} + q
\sigma_{(1,1)} + $

\hspace{2.3cm}$(T_1 - T_4)q \sigma_{(1)} +
 (T_1 - T_2)(T_3- T_4)q$

$ \sigma_{(2,1)} \circ \sigma_{(2,2)} =
(T_1-T_4)(T_1-T_3)(T_2-T_4)\sigma_{(2,2)} + q\sigma_{(2,1)} +
(T_1-T_3)q \sigma_{(2)}+$

\hspace{2.3cm}  $(T_2-T_4)q \sigma_{(1,1)} + (T_1 - T_3)(T_2-T_4)q
\sigma_{(1)}$

$\sigma_{(2,2)} \circ \sigma_{(2,2)} =
(T_1-T_4)(T_1-T_3)(T_2-T_4)(T_2-T_3)\sigma_{(2,2)} +
(T_1+T_2-T_3-T_4)q\sigma_{(2,1)} +$

\hspace{2.3cm} $(T_1-T_3)(T_2-T_3)q \sigma_{(2)}+
(T_2-T_3)(T_2-T_4)q \sigma_{(1,1)} + $

\hspace{2.3cm} $(T_1 - T_3)(T_2-T_3)(T_2-T_4)q \sigma_{(1)} + q^2$

\subsection{Final remarks:} 1. {\em Presentation and equivariant
quantum Giambelli.} A presentation for the equivariant quantum
cohomology of partial flag manifolds was obtained in \cite{Kim1},
Thm. I. One would like to find polynomial representatives
corresponding to the equivariant quantum Schubert classes
$\sigma_\lambda$ (i.e. an equivariant quantum Giambelli formula).
In a paper in preparation (\cite{Mi1}, we obtain another
presentation for the equivariant quantum cohomology of the
Grassmannians, and we find such an equivariant quantum Giambelli
formula.

2. {\em Positivity.} In a paper in preparation (\cite{Mi2}), we
prove that the EQLR coefficients {\it for any homogeneous space
$G/P$} enjoy the same positivity property as the equivariant
coefficients (see \cite{Gr} for the latter). This implies that the
EQLR coefficients from this paper are homogeneous polynomials in
variables $T_1-T_2,..., T_{m-1}-T_m$ with \textit{nonnegative}
coefficients.

3. {\em Generalization to other homogeneous spaces.} The fact that
EQ Pieri/Monk rule doesn't contain any mixed terms generalizes to
all homogeneous spaces $G/P$ (see remark after Cor. \ref{Pieri
vanishing}). The proof of this result will be given elsewhere.

\section{Appendix} The Appendix contains the definition and some
properties of the equivariant Gysin maps used in \S 2.2. It also
includes a proof of Proposition \ref{eqqcoh}.

\subsection{Equivariant Gysin morphisms} Let $f:X \longrightarrow Y$ be a
morphism of projective varieties, with $Y$ smooth. Let $d= \dim
(X) - \dim (Y) $ (complex dimensions). Define a Gysin map
$f_\star:H^i(X) \longrightarrow H^{i - 2d} (Y)$ by the composite
$$ \begin{CD} H^i(X) @> { \cap [X] } >> H_{2 \dim (X) - i}
(X) @> {f_\star} >> H_{ 2 \dim (X) - i} (Y) \simeq H^{i - 2d}(Y)
\end{CD} $$ where $[X]$ is the fundamental class of $X$ in the singular homology
group $H_{2 \dim X}(X)$, and the middle $f_\star$ is the singular
homology push-forward (if $X$ or $Y$ were not compact, one should
use Borel-Moore homology). The last isomorphism is given by
Poincar\'e duality. We need the following property of the Gysin
map:

\begin{lemma}\label{fiber square} Consider the following fiber square of projective
varieties:
$$ \begin{CD} X' @>{i}>> X \\ @V {f'}VV @V {f} VV \\ Y' @>{j}>>Y
\end{CD} $$ where $Y,Y'$ are smooth and $i,j$ are regular
embeddings of the same (complex) codimension $c$. Then $f'_\star
i^\star = j^\star f_\star$ as maps $H^i(X) \to H^{i-2d}(Y')$.
\end{lemma}

\begin{proof} The proof is given in my thesis \cite{Mi} (one could also see \cite{FM}).
\end{proof}

Assume the map $f:X \to Y$ (with $Y$ smooth) is $T-$equivariant.
Then it determines a Gysin map of the cohomology of the
finite-dimensional approximations $f_{\star,n}:H^i(X_{T,n})
\longrightarrow H^{i-2d} (Y_{T,n})$. Define the equivariant Gysin
map $f_\star^T:H^\star_T(X) \longrightarrow H^{i-2d}_T(Y)$ as the
unique map that makes the following diagram commute:
$$ \begin{CD} H^i(X_{T,n}) @<{res}<<H^i_T(X) \\
@V{f_{\star,n}}VV @V{f_\star^T}VV\\ H^{i- 2d}(Y_{T,n})
@<{res}<<H^{i-2d}_T(Y) \end{CD} $$for any integer $n$. The
horizontal maps $res$ are the cohomology pull-backs induced by the
inclusions $X_{T,n} \to X_T$ (resp. $Y_{T,n} \to Y_T$). The
uniqueness of $f_\star^T$ follows from the fact that the
equivariant cohomology can be computed by passing to the limit on
the ordinary cohomology of the finite dimensional approximations
(see \S 2.2). The fact that such a definition is independent of
the choice of the approximation $f_n:X_{T,n} \longrightarrow
Y_{T,n}$ follows by applying Lemma \ref{fiber square} to the fiber
square \begin{equation*}
\begin{CD}
X_{T, n_1} @> {i}>>X_{T, n_2}\\@V{f_{n_1}}VV @V{f_{n_2}}VV \\
Y_{T, n_1} @>{j}>> Y_{T, n_2} \end{CD} \end{equation*} for
integers $n_1 < n_2$.

Another property of the equivariant Gysin map is its compatibility
with the restriction to the fiber:
\begin{equation}\label{compatibility} \begin{CD} H^i(X) @<{res}<<H^i_T(X) \\
@V{f_{\star}}VV @V{f_\star^T}VV\\ H^{i- 2d}(Y)
@<{res}<<H^{i-2d}_T(Y)
\end{CD} \end{equation} This follows by applying again Lemma \ref{fiber square} to
the fiber square $$ \begin{CD}
X @> {i}>>X_{T, n}\\ @V{f}VV @V{f_{n}}VV \\
Y @>{j}>> Y_{T, n} \end{CD} $$

\begin{proof}[\bf{Proof of Prop. \ref{eqqcoh}}]
Part (a) of the proposition is proved in \cite{Kim2}, \S 3.3 , using a
slightly different definition of the EQLR
coefficients.\begin{footnote}[6]{The main part of Kim's proof is the
associativity of the equivariant quantum cohomology. This follows from
the proof of the associativity of the (non-equivariant) quantum
cohomology, by noting that all the maps involved in the proof are in fact
$T-$equivariant. For details see \cite{Kim2}, \S 3.3.}\end{footnote} For
the sake of completeness we recall Kim's definition, and we prove it is
equivalent to ours.

Let $a_1,a_2,a_3$ be three equivariant cohomology classes in
$H^\star_T(X)$. The equivariant Gromov-Witten invariant, denoted
$I_{3,d}^{X_T}(a_1,a_2,a_3)$ (Kim's notation), is defined by
\[I_{3,d}^{X_T}(a_1,a_2,a_3)=\pi_\star^T((ev_1^T)^\star(a_1) \cup
(ev_2^T)^\star(a_2) \cup (ev_3^T)^\star(a_3)) \] where
$\pi_\star^T$ is the equivariant Gysin morphism (\S 2.2). Let $ <
\cdot , \cdot
> $ be the $q-$linear extension of the equivariant Poincar\'e
pairing defined in section 2.2. Then, according to \cite{Kim2}, \S
3.3 , the equivariant quantum multiplication is the unique
multiplication (denoted $\diamond$) such that
\[ <a_1 \diamond a_2, a_3> = \sum_d q^d I_{3,d}^{X_T}(a_1,a_2,a_3) \]
By taking $a_1:=\sigma_\lambda^T$, $a_2:=\sigma_\mu^T$, Kim's
definition implies that the coefficient of $\sigma_\nu^T$ in
$\sigma_\lambda^T \diamond \sigma_\mu^T$ is equal to
\[\sum_d q^d
I_{3,d}^{X_T}(\sigma_\lambda^T,\sigma_\mu^T,(\sigma_\nu^T)^\vee)\]
where $(\sigma_\nu^T)^\vee$ is the dual of $\sigma_\nu^T$ with
respect to the equivariant Poincar\'e pairing. But equivariant
duality (Proposition 2.2), implies that $(\sigma_\nu^T)^\vee =
\widetilde{\sigma}_{\nu^\vee}^T$, which shows that the
multiplications $\diamond$ and $\circ$ coincide.

Part (b) is a consequence of the fact that the EQLR coefficients
specialize to both equivariant and quantum ones. For
the convenience of the reader, we sketch a proof of this fact.\\

\noindent \textit{Claim 1.} The EQLR coefficient
$c_{\lambda,\mu}^{\nu,0}$ is equal to the equivariant LR
coefficient $c_{\lambda,\mu}^{\nu}$.
\begin{proof} If $d=0$, $\overline{\mathcal{M}}_{0,3}(X,0)_T = X_T$,
so the definition of the EQLR becomes the definition of the
equivariant LR coefficients from the end of Section 2.2.
\end{proof}

\noindent \textit{Claim 2.} If $|\lambda|+|\mu|=|\nu|+md$, the
EQLR coefficient $c_{\lambda,\mu}^{\nu,d}$ is equal to the quantum
LR-coefficient.

\noindent {\it Proof.} Denote by $D= p(m-p)+md$ the dimension of
$\mbotreix$. The hypothesis implies that
$|\lambda|+|\mu|+|\nu^\vee|=D$. We use the compatibility between
the ordinary and equivariant Gysin maps (see diagram
(\ref{compatibility}):
$$
\begin{CD}
 H^{2D}(\mbotreix) @<{res}<< H^{2D}_T(\mbotreix) \\ @ V{\pi_\star} VV @V {\pi_\star^T}VV \\
 H^{0} (pt) @<
 {res}<< H^{0}_T(pt)
\end{CD} $$ The top restriction map sends the classes $ev_i^{T\star}(\sigma_\lambda^T)$
and $ev_i^{T\star}(\widetilde{\sigma}_\lambda^T)$ to
$ev_i^\star(\sigma_\lambda)$, ($i=1,2,3$) and the bottom
restriction map is an isomorphism. The claim now follows from
Prop. \ref{qdef}.
\end{proof}

\noindent \textit{Remark:} Another proof of Claim 2, using the
``restriction'' property  of the equivariant quantum cohomology,
can be found in \cite{Kim1} Prop. 1, \S 5 (see also \S 4.4).


\begin{thebibliography}{[AS]}

\bibitem [AB] {AB} Atyiah, M., F., Bott, R., The moment map and
equivariant cohomology, {\em Topology} {\bf23} (1984) no. 1, 1- 28

\bibitem [AS]{AS} Astashkhevich, A., Sadov, V., Quantum cohomology
of partial flag manifolds $F_{n_1,...,n_k}$, {\em Commun. Math.
Phys.} {\bf 170} (1995), 503-528

\bibitem[Be] {Be} Bertram, A., Quantum Schubert calculus, {\em Adv.
Math}, {\bf 128} (1997), no. 2, 289-305

\bibitem[BCF] {BCF}Bertram, A., Ciocan-Fontanine, I., Fulton, W., Quantum
multiplication of Schur polynomials, {\em Journal of Algebra} {\bf
219} (1999), no.2, 728-746

\bibitem[Br1]{Br1} Brion, M., Equivariant cohomology and
Equivariant Intersection Theory  {\em NATO Adv. Sci. Inst. Ser. C
Math. Phys. Sci.}, {\bf 514}, Representation theories and
algebraic geometry (Montreal, PQ, 1997), 1--37, Kluwer Acad.
Publ., Dordrecht, 1998.

\bibitem[Br2] {Br} Brion, M., Poincar\'e Duality and Equivariant
(Co)homology, {\em Michigan Math. J. - Special volume in honor of
William Fulton} {\bf 48} 2000, 77-92

\bibitem [Bu1] {Bu1} Buch, A. S., Quantum cohomology of Grassmannians, {\em to appear
in Compositio Math.}, arXiv: math.AG/0106268

\bibitem[BKT]{BKT1} Buch, A. S., Kresch A., Tamvakis H., Gromov-Witten
invariants on Grassmannians, {\em to appear on J. Amer. Math. Soc.},
arXiv: math.AG/0306388

\bibitem[C]{C} Ciocan-Fontanine, I., Quantum cohomology of flag varieties,
{\em Internat. Math. Res. Notices} (1995), 263-277

\bibitem[EG]{EG} Edidin, D., Graham, W., Equivariant intersection
Theory (with an appendix by Angelo Vistoli: The Chow ring of
$\mathcal{M}_2$), {\em Invent. Math.}{\bf 131}, (1998) 595-634

\bibitem[FGP]{FGP} Fomin, S., Gelfand, S., Postnikov, A., Quantum Schubert
Polynomials, {\em J. Amer. Math. Soc.} {\bf 10} (1997), 565-596

\bibitem[F1]{F1}Fulton, W., Young Tableaux, {\em Cambridge University
Press}, Cambridge, 1997

\bibitem[F2]{F3} Fulton, W., Intersection Theory, {\em Springer
Verlag} 2nd edition (1998)

\bibitem[FP]{FP} Fulton, W., Pandharipande, R., Notes On Stable Maps And
Quantum Cohomology, {\em Proc. Sympos. Pure Math.} {\bf 62}, Part
2, Amer. Math.  Soc., Providence, RI, 1997.

\bibitem[FM] {FM} Fulton, W., MacPherson, R., Categorical
framework for the study of singular spaces, {\em Memoirs of AMS},
{\bf 1981} vol. 31 no. 243

\bibitem[FW] {FW} Fulton, W., Woodward, C., On the quantum product
of Schubert classes, {\em to appear in J. of Alg. Geom.}, arXiv:
math.AG/0112183

\bibitem[G]{G} Givental, A., Equivariant Gromov-Witten invariants,
{\em IMRN} (1996) 613-663

\bibitem[GK]{GK} Givental, A., Kim, B., Quantum cohomology of flag
manifolds and Toda lattices, {\em Comm. Math. Phys.} {\bf 168}
(1995), 609-641

\bibitem[GKM]{GKM} Goreski, M., Kottwitz, R., MacPherson, R.,
Equivariant cohomology, Koszul duality, and the localization
theorem, {\em Invent. Math.} {\bf 131} (1998), no. 1, 25-83

\bibitem[Gr]{Gr} Graham, W., Positivity in equivariant Schubert
calculus, {\em Duke Math. J.} {\bf 109} (2001), no. 3, 599-614

\bibitem[H]{H} Husemoller, D., Fibre Bundles, {\em Springer-Verlag}
(1975) 2nd edition

\bibitem[Kim1]{Kim1}Kim, B., Quantum Cohomology of partial flag
manifolds and a residue formula for their intersection pairings,
{\em IMRN} {\bf 1995}, no.1, 1-15

\bibitem[Kim2]{Kim2} Kim, B.,On equivariant quantum cohomology,
{\em IMRN} {\bf 17} (1996), 841-851

\bibitem[Kim3]{Kim3} Kim, B., Quantum cohomology of flag manifolds {G/B}
and quantum Toda lattices, {\em Annals of Math.} {\bf 149} (1999),
129-148

\bibitem[KiMa]{KiMa} Kirillov, A. N., Maeno, T., Quantum double
Schubert polynomials, quantum Schubert polynomials and
Vafa-Intriligator formula, {\em Discrete Math.} {\bf 217} (2000)
no. 1-3, 191-223

\bibitem[KM]{KM} Kontsevich,M., Manin, Y., Gromov-Witten classes,
quantum cohomology and enumerative geometry, {\em Comm. Math.
Phys.} {\bf 164} (1994), 525-562

\bibitem[KT]{KT1} Knutson, A., Tao, T., Puzzles and equivariant
cohomology of Grassmannians, {\em Duke Math. J.} {\bf 119} (2003)
issue 2, 221-260

\bibitem[Mi]{Mi} Mihalcea, L.C., Ph.D. Thesis, University of
Michigan

\bibitem[Mi1] {Mi1} Mihalcea, L. C., Polynomial representatives for the Schubert classes in the
equivariant (quantum) cohomology  of the Grassmannian {\it in
preparation}

\bibitem[Mi2]{Mi2} Mihalcea, L. C., Positivity in equivariant quantum
Schubert calculus {\em in preparation}

\bibitem[MS]{MS} Molev, A. I., Sagan B., A Littlewood-Richardson rule
for factorial Schur functions, {\em Trans. Amer. Math. Soc.} {\bf
351} (1999), no. 11, 4429-4443

\bibitem[O]{O} Okounkov A., Quantum immanants and higher Capelli
identities, {\em Transformations Groups} {\bf 1} (1996), 99-126

\bibitem[Po]{Po} Postnikov, A., Affine approach to quantum Schubert calculus,
{\em preprint arXiv: math.CO/0205165}


\bibitem[QC]{qcoh} Quantum Cohomology at the Mittag-Leffler Institute,
{\em edited by Paolo Aluffi}, 1996


\bibitem[R]{R} Robinson, S., A Pieri-type formula for the
equivariant cohomology of the flag manifold, {\em Journal of
Algebra} {\bf 249} , 38-58 (2002).

\bibitem[S]{S} Stembridge, J.R., A concise proof of the
Littlewood-Richardson rule, {\em Electron. J. Comb.} {\bf 9} 2002

\bibitem[Y]{Y} Yong, A., Degree bounds in Schubert calculus, {\em
Proceedings of the AMS}, Vol. 131, Number 9, 2649-2655 (2003).

\bibitem[W]{W} Witten, E., The Verlinde algebra and the cohomology of the
Grassmannian, {\em Geometry, Topology and Physics}, Internat.
Press, Cambridge, MA, 1995, 357-422
\end{thebibliography}
\end{document}